\documentclass[11pt,leqno]{article}
\usepackage{amsthm,amsmath,amsfonts,latexsym,amssymb}
\usepackage{mathrsfs,dutchcal,arydshln,gensymb,appendix,abstract}
\usepackage[dvips]{graphics,color}
\usepackage{graphicx}
\usepackage{tikz}
\date{}

\begin{document}
	\title{On FKM isoparametric hypersurfaces in $\mathbb{S}^n \times \mathbb{S}^n$ and new area-minimizing cones}  
	\author{\centerline {Hongbin Cui}}
 
 \maketitle

\begin{abstract}
We present two generalizations for the celebrated works of Ferus-Karcher-Münzner \cite{FKM81} and Wang \cite{W94}. We first show that an isoparametric foliation on $\mathbb{S}^{2n+1}$ constructed by Ferus-Karcher-Münzner naturally yields an isoparametric foliation on its submanifold $\mathbb{S}^n \times \mathbb{S}^n$ with one same focal variety. The second part concerns area-minimizing cones; all known regular area-minimizing hypercones are realized as real algebraic varieties: isoparametric cones (cf. \cite{W94}). As a noteworthy application, we extend area-minimizing isoparametric hypercones in \cite{W94} to codimension-two cases, and obtain infinitely many families (each containing infinitely many members) of area-minimizing subcones of Simons cones. 
\end{abstract}

		
		
		\bigskip
		
		\tableofcontents

\section{Introduction}

\medskip\noindent

The generalized Plateau's problem stands as a classic topic in geometric measure theory that seeks to find an $k$-dimensional area-minimizing surface (or integral currents, cf. \cite{FF60}) $M$ with a given $(k-1)$-dimensional boundary $\Sigma$ in some Euclidean space $\mathbb{R}^{l+1}$. Singularities do appear in very natural settings for area-minimizing surfaces, whether the boundary has singularities or the ambient space is higher-dimensional. In \cite{DeL23}, De Lellis gave a rather extensive review of Plateau's problem, its recent developments, and some of the remaining challenges. Several programmatic problems that are central to the regularity theory are \cite[p. 875]{DeL23}:

\begin{itemize}
\item \textbf{Q1:} Under which conditions singularities can be ruled out?
\item \textbf{Q2:} How large can the set of singularities be?
\item \textbf{Q3:} Which structural properties can the singularities have? 
\end{itemize}

 Towards \textbf{Q2} for interior singularities,  the celebrated interior regularity theorems provide quantitative descriptions for the Hausdorff dimension of the interior singular set (see a nice survey in \cite{M16}), the object of core significance for the codimension $c:=l+1-k=1$ cases is the following Simons cones, cf. \cite{Fe69}, \cite{Fe70}. Meanwhile, as in \textbf{Q3}, understanding the concrete structures of the singularities is of equal importance. Many significant advances on \textbf{Q3} and  \textbf{Q2} were also reviewed in \cite{DeL23}. More recently, a striking work of Liu \cite{L25a} shows that there exist area-minimizing surfaces with fractional-dimensional singular sets. Another way to think about the shape of these singularities may be their first-order behaviour: \textit{tangent cones}, as a generalization of the tangent plane at regular points. The tangent cones of area-minimizing surfaces at their interior singularities are also area-minimizing (cf. \cite[Theorem 5.4.3]{Fe69}). 

\textbf{Example 1} (Simons cones). The cones
$$
\mathcal{C}_{n}:=C(\mathbb{S}^n \times \mathbb{S}^n)=\{(x,y)\in \mathbb{R}^{n+1}\times \mathbb{R}^{n+1}: |x|=|y| \}, \quad \ n\geq 3
$$ 
are area-minimizing in $\mathbb{R}^{2n+2}$ with an isolated singularity at the origin. $\mathcal{C}_{n}$ are the first family known area-minimizing hypercones which were proven by Bombieri, De Giorgi and Giusti in \cite{BDGG69} (also see \cite{La72}, \cite{Fe74}, \cite{La98}, \cite{Da04}, \cite{DP09} for alternative proofs). Together with the previous works in \cite{F62}, \cite{G65}, \cite{Alm66}, \cite{S68}, it gives a complete answer to the famous Bernstein problem. We note here, the intersections of $\mathcal{C}_{n}$ with $\mathbb{S}^{2n+1}(\sqrt{2})$, called the \textit{links}, are the products of spheres
$$
\mathbb{S}^{n}(1) \times \mathbb{S}^{n}(1) \subset \mathbb{S}^{2n+1}(\sqrt{2}).
$$ 

Towards \textbf{Q1}, Brian White \cite{B86} suggested that finding the smooth area-minimizing surfaces under the generic metric, i.e., a perturbation of the ambient metric. In a recent study, Liu \cite{L25b} has made remarkable contributions to the study of \textbf{Q1} in general cases. He proved that when $k\geq5, c\geq 3$ or $3\leq c \leq k \leq 4$, the singularities cannot be perturbed under generic metrics.  Liu's research highlights the differences between the cases of \textit{low codimension} $c\leq 2$ and higher codimensions, and several area-minimizing cones of codimension $c=3$ play key roles in his study, cf. \cite[fact 1.1.2 and fact. 1.1.3]{L25b}. In general dimensional cases, Liu conjectured that area-minimizing surfaces of codimension $c=1$ or $c=2$ are smooth in generic metrics. For this purpose, examples of regular area-minimizing cones of low dimensional $c\leq 2$ are clearly precious. We refer the reader to the Appendix for more information on  regular area-minimizing cones with low codimension that are not complex varieties.

\medskip\noindent

In this paper, we construct infinitely many families (each containing infinitely many members) of non-holomorphic area-minimizing cones with codimension $c=2$ (cf. Theorem B) in $\mathbb{R}^{2n+2}$ for $n\geq 11$. They are the subcones of Simons cones: the intersections of Simons cones with some specific FKM isoparametric cones. Our results are based on an natural generalization for the celebrated work of Ferus, Karcher and Münzner \cite{FKM81} on isoparametric hypersurfaces in the spheres, see Theorem A below. 

\medskip\noindent

\subsection{The construction of FKM isoparametric hypersurfaces in $\mathbb{S}^{n} \times \mathbb{S}^{n}$ and the minimality}

\medskip\noindent

We focus on the minimal isoparametric hypersurfaces in the links of Simons cones: $\mathbb{S}^n \times \mathbb{S}^n$. A key feature of $\mathbb{S}^n \times \mathbb{S}^n$ is its canonical product involution: $P(v_1,v_2)=(v_1,-v_2)$. In the context of $P$, for hypersurfaces in $\mathbb{S}^n \times \mathbb{S}^n$ with a unit normal field $N$, F. Urbano \cite{Ur19} defined $C:=\langle P N,N \rangle$ which is independent of the choice of orientation. The following simple lemma illustrates the relationship between the product angle function $C$ and the simultaneous minimality property of hypersurfaces in $\mathbb{S}^{n}\left(1\right) \times \mathbb{S}^{n}\left(1\right) \subset \mathbb{S}^{2n+1}(\sqrt{2})$.

\medskip\noindent

\textbf{Lemma 1.1 (See \cite[Lemma 2.4]{Chen21})} \textit{Let $M$ be a hypersurface of $\mathbb{S}^{n}\left(1\right) \times \mathbb{S}^{n}\left(1\right) \subset \mathbb{R}^{2n+2}$, let $H$ and $\widetilde{H}$ denote the mean curvature vector of $M$ in $\mathbb{S}^{n}\left(1\right) \times \mathbb{S}^{n}\left(1\right)$ and $\mathbb{S}^{2n+1}(\sqrt{2})$, $(x,y)$ denote the position vector. Then 
$$
\widetilde{H}(x,y)=H(x,y)+\frac{C}{2}(x,-y).
 $$}

\medskip\noindent

Thus, a hypersurface $M \subset \mathbb{S}^{n} \times \mathbb{S}^{n}$ is minimal in $\mathbb{S}^{2n+1}(\sqrt{2})$ if and only if $M$ is minimal in $\mathbb{S}^{n} \times \mathbb{S}^{n}$ with $C\equiv0$. More generally, if $C$ is constant, then $M$ has constant mean curvature (CMC) in $\mathbb{S}^{2n+1}(\sqrt{2})$ if and only if $M$ has CMC in $\mathbb{S}^{n} \times \mathbb{S}^{n}$. Minimal hypersurfaces in $\mathbb{S}^{2} \times \mathbb{S}^{2}$ satisfying $C\equiv0$ were classified in \cite[Theorem 3.5]{LVWY24}.

\medskip\noindent

For isoparametric hypresurfaces in $\mathbb{S}^{n} \times \mathbb{S}^{n}$, Urbano \cite{Ur19} constructed a family of examples in $\mathbb{S}^{n} \times \mathbb{S}^{n}$:
$M_{t}=\{(x,y)\in \mathbb{S}^n \times \mathbb{S}^n| \langle x,y \rangle=t\}$. In fact, priori to the research of Urbano, Qian-Tang (\cite[Theorem 1]{QT16}) had already obtained these examples among their many results. Their examples are homogeneous, it is clear that $SO(n+1)$ acts transitively by isometries on $M_t$ by
$
A(p, q)=(A p, A q), \quad A \in S O(n+1),
$
and hence $\left\{M_t, t \in(-1,1)\right\}$ is a family of homogeneous hypersurfaces. In fact, $\mathbb{S}^{n} \times \mathbb{S}^{n}$ admints an isometric group $(O(n+1) \times O(n+1)) \rtimes \mathbb{Z}_2$ where $\mathbb{Z}_2$ action is given by $\tau (p,q)=(q,p)$. These are the only known examples of isoparametric hypersurfaces in $\mathbb{S}^{n} \times \mathbb{S}^{n}$ before. 

\medskip\noindent

We now give a far-reaching generalization of Qian-Tang and Urbano's isoparametric examples from the further restrictions of FKM isoparametric hypersurfaces on $\mathbb{S}^{n}(1) \times \mathbb{S}^{n}(1) \subset \mathbb{S}^{2n+1}(\sqrt{2})$. In \cite{FKM81}, Ferus, Karcher, and Münzner constructed new examples of isoparametric hypersurfaces in spheres $\mathbb{S}^{2n+1}$  with $4$ distinct principal curvatures using a (symmetric) Clifford system $P_0, \cdots, P_{m}$ in $\mathbb{R}^{2n+2}$. They define the polynomials as 
$$
\mathbf{F}(\widetilde{x})= \langle \widetilde{x}, \widetilde{x} \rangle^2-2 \sum_{i=0}^{m} \langle P_i \widetilde{x}, \widetilde{x} \rangle^2, \ \ \ m\geq 1
$$
which we refer to as \textit{FKM isoparametric polynomial} and their regular level sets as \textit{FKM isoparametric hypersurfaces}. There are infinitely many families, each with infinitely many members, of non-homogeneous isoparametric hypersurfaces in their examples (cf. \cite[p. 177-178]{CR15}).  The FKM  isoparametric hypersurfaces are of remarkable importance, since after decades of continuous work by many mathematicians, a complete classification of isoparametric hypersurfaces in spheres was finally accomplished in recent years (cf. \cite{CCJ07}, \cite{Chi20} or \cite{Chi18}): they are either homogeneous or FKM type.

\medskip\noindent

By means of the eigenvalue space decomposition of $P_{0}$, we can set $\widetilde{x}=(x,y) \in \mathbb{R}^{n+1} \times \mathbb{R}^{n+1}$, and 
$$
P_0=\left(\begin{array}{rr}
I_l & 0 \\
0 & -I_l
\end{array}\right), \quad P_1=\left(\begin{array}{cc}
0 & I_l \\
I_l & 0
\end{array}\right), \quad P_{1+\alpha}=\left(\begin{array}{rr}
0 & A_\alpha \\
-A_\alpha & 0
\end{array}\right) \quad \text { for } 1 \leq \alpha \leq m-1,
$$
where $\left\{A_1, A_2, \ldots, A_{m-1}\right\}$ satisfy
\begin{equation}\label{skc0}
    \left\{\begin{array}{c}
A_p A_q+A_q A_p=-2 \delta_{p q} I \\
A_p=-A_p^{T}
\end{array}\right.
\end{equation}
are skew-symmetric representations of the Clifford algebra $Cl_{m-1}$ on $\mathbb{R}^{n+1}$, see Sect. 2.2 for more details. Then
\begin{equation}\notag
    \mathbf{F}(\widetilde{x})=(|x|^2+|y|^2)^2-2(|x|^2-|y|^2)^2-8\langle x, y\rangle^2-8\sum_{\alpha=1}^{m-1}\langle A_\alpha x, y\rangle^2,
\end{equation}
the restriction of $\mathbf{F}$ on $\mathbb{S}^{n} \times \mathbb{S}^{n} \subset \mathbb{S}^{2n+1} $ is equivalent to the restriction of the following bihomogeneous polynomial:
$$
F(x,y)=\langle  x, y\rangle^2+\sum\limits_{q=1}^{m-1}\langle A_q x, y\rangle^2,
$$
We refer to this as the splitting form of the FKM isoparametric polynomials. We then make precise that 

\medskip\noindent

\textbf{Theorem A} \textit{
The level sets
$$
 M_{t}=\{ (x,y)\in \mathbb{S}^n(1) \times \mathbb{S}^n(1)| \langle  x, y\rangle^2+\sum\limits_{q=1}^{m-1}\langle A_q x, y\rangle^2=t \}, t\in (0,1)
$$
gives isoparametric hypersurfaces in $\mathbb{S}^n(1) \times \mathbb{S}^n(1)$ satisfying $C\equiv 0$,
where $\left\{A_1, A_2, \ldots, A_{m-1}\right\}$ are skew-symmetric representations of the Clifford algebra $Cl_{m-1}$ on $\mathbb{R}^{n+1}$.}

\medskip\noindent

For $m=1$, it recovers Qian-Tang and Urbano's examples which can be defined on any $\mathbb{S}^n \times \mathbb{S}^n$. For $m\geq 2$,  they defined on $\mathbb{S}^{2k-1} \times \mathbb{S}^{2k-1}$ which depends on the dimension of the irreducible representation of Clifford algebra (see Sect. 2.2).  Clearly, these examples can be defined in $\mathbb{S}^n(a) \times \mathbb{S}^n(b)$ for $t \in (0, a^2b^2)$ by virtue of their bihomogeneous property (see Proposition 2.8). 

\medskip\noindent

\textbf{Definition 1.2.} \textit{We call those new isoparametric hypersurfaces in Theorem A  with $m\geq 2$, the FKM isoparametric hypersurfaces in $\mathbb{S}^{n} \times \mathbb{S}^{n}$.}

\medskip\noindent

\textbf{Remark 1.3.} \textit{There has been a growing interest in the classification problem for isoparametric hypersurfaces in the products of space forms (especially in terms of establishing the constant angle properties, for the newest progress, see \cite{DP24}, \cite{TXY25}, \cite{DP25}). We remark that, except for the new examples presented here, further new advances in the classification of isoparametric hypersurfaces in $\mathbb{S}^n \times \mathbb{S}^n$ have been obtained in more recent independent works, cf. \cite{TXY26+} and \cite{W26+}.}

\medskip\noindent
\medskip\noindent

\subsection{New area-minimizing cones and discussions}

\medskip\noindent

$\mathbb{S}^n \times \mathbb{S}^n \ (n\geq 2)$ admits positive Ricci curvature, by a well-known analysis for the associated Riccati type equation as in \cite{G04} or \cite{GT14}, there is a unique minimal isoparametric hypersurface $\Sigma$ within our isoparametric family. By Theorem A and Lemma 1.1, $\Sigma$ is also minimal in $\mathbb{S}^{2n+1}$.

\medskip\noindent

For $m=1$, it is given by Qian-Tang and Urbano's example $M_{0}=\{(x,y)\in \mathbb{S}^n \times \mathbb{S}^n| \langle x,y \rangle=0\}$, which can be identified as a Stiefel manifold of orthogonal 2-frames. Let $z_i=x_i+\sqrt{-1} y_i \ (1\leq i \leq n+1)$, then $z_{i}^2=x_{i}^2-y_{i}^2+2\sqrt{-1}x_{i}y_{i}$. It follows that the cones over their examples
$$
\mathcal{C}(M_{0})=
\left\{\left(z_1, \cdots, z_{n+1}\right) \in \mathbb{C}^{n+1} \mid \sum_{i=1}^{n+1} z_i{ }^2=0\right\}
$$
are holomorphic varieties (complex hyperquadrics) that singular only at the origin, thus area-minimizing (cf. \cite{HL82}). 

\medskip\noindent

For the general cases $m\geq 2$, the minimal FKM isoparametric hypersurfaces are given by 
$$
\Sigma=F^{-1} \left(\frac{m-1}{n-1}\right) \cap \mathbb{S}^n \times \mathbb{S}^n,
$$
thus $\mathcal{C}(\Sigma)$ are minimal cones in $\mathbb{R}^{2n+2}$ which are contained in the Simons cones $\mathcal{C}_{k\delta(m)-1}=\mathcal{C}(\mathbb{S}^{k\delta(m)-1} \times \mathbb{S}^{k\delta(m)-1})$. In the second part of this paper, we establish 

\medskip\noindent

\textbf{Theorem B.} \textit{The minimal cones $\mathcal{C}(\Sigma)\subset \mathcal{C}_{k\delta(m)-1}$
$$
\mathcal{C}(\Sigma)=\left\{(x,y)\in \mathbb{R}^{k\delta(m)} \times \mathbb{R}^{k\delta(m)}:   \langle x, y\rangle^2+\sum_{q=1}^{m-1}\langle A_q x, y\rangle^2=\frac{m-1}{n-1}|x|^2|y|^2 \ \ {\rm and} \ \ |x|=|y| \right\},
$$ 
are area-minimizing varieties if (i): in the irreducible case $k=1$, $ \delta(m)\geq 16$; (ii): in the reducible case $k\geq 2$, $ k\delta(m)\geq 12$. Where $\left\{A_1, A_2, \ldots, A_{m-1}\right\}$ are skew-symmetric representations of the Clifford algebra which satisfying \eqref{skc0}, the values of $\delta(m)$ are given in \eqref{table}.}

\medskip\noindent

The proof of Theorem B relies on Lawlor's curvature criterion (see Sect. 3.1). Specifically, we need to calculate two quantities of $\mathcal{C}(\Sigma)$: the vanishing angles (see Def. 3.4) and the normal radius (see Def. 3.1). Then if two times of vanishing angles are less than the normal radii, we can prove that $\mathcal{C}(\Sigma)$ is area-minimizing (see Theorem 3.5). The computations of the normal radius are carried out in Sect. 3.3 by repeatedly applying \eqref{skc0} to find the shortest normal geodesic of $\Sigma$ in $\mathbb{S}^{2n+1}(\sqrt{2})$. 

\medskip\noindent

\textbf{The second fundamental forms}. To estimate the vanishing angles, a key challenge is to determine the maximum value $\alpha^2$ of the squared norm of the shape operator for $\frac{1}{\sqrt{2}}\Sigma \subset \mathbb{S}^{2n+1}(1)$. First, by applying a computation through Bochner formula (cf. \cite[Prop. 1.47]{CM11}), we directly derive $|B|^2$ for the squared norm of the shape operator for $\frac{1}{\sqrt{2}}\Sigma \subset \mathbb{S}^{n}(\frac{1}{\sqrt{2}}) \times \mathbb{S}^{n}(\frac{1}{\sqrt{2}})$,
which avoids the complex analysis of principal curvatures employed in \cite{Ur19}\footnote{In fact, as Urbano's computation for his isoparametric examples in $\mathbb{S}^2 \times \mathbb{S}^2$, the principal direction is found based on the special complex structures and it cann't split properly \cite{DP25}, thus determining the concrete values for principle curvatures for those more general FKM isoparametric hypersurfaces are not easy. There are some cross terms as in \eqref{crt}.}. We then need to establish that $\alpha^2=|B|^2$ (See Sect. 3.2). To tackle this problem, we slice the hypersurface $\Sigma \subset \mathbb{S}^{n} \times \mathbb{S}^{n}$ starting from a given point $(x,y)\in \Sigma$:
$$\Sigma_{1}: \Sigma \cap\left(S^{n} \times\{y\}\right) \hookrightarrow S^{n} \times\{y\} \ \ {\rm and} \ \ \Sigma_{2}: \left(\{x\}\times S^{n}\right) \cap \Sigma \hookrightarrow\{x\} \times S^{n}.$$
we prove that both $\Sigma_1$ and $\Sigma_{2}$ are \textit{minimal isoparametric hypersurfaces} in $\mathbb{S}^{n}$. We then investigate the relationship between the second fundamental forms of $\Sigma$ and those of $\Sigma_1$ and $\Sigma_2$: apart from the tangent vectors to the two slices of $\Sigma$, there exists only one additional tangent vector $T$ for $\Sigma$. By Urbano's formula $\bar\nabla C=-2A(T)$ \cite[Lemma 1(1)]{Ur19}, all non-vanishing components of the second fundamental forms
belong to the two slices. These non-vanishing terms can be computed by considering curves that move along one single slice. Throughout these analyses, we ultimately establish that $\alpha^2=|B|^2$.

\medskip\noindent

In Sect. 4, inspired by the construction of minimal product cones in \cite{TZ20}, we establish that a $k$-dimensional minimal product cone over FKM isoparametric hypersurfaces in products of two spheres is necessarily area-minimizing when $k\geq 21$.

\medskip\noindent

\textbf{Organizations of the paper.} In Sect. 2.1, we investigate the geometry of level set hypersurfaces in $\mathbb{S}^n \times \mathbb{S}^n$. The most important part is that we define a family of canonical frames (see Definition 2.1) for those hypersurfaces, then we can derive a formula for the squared length of the second fundamental form via the Bochner formula and these canonical frames. In Sect. 2.2, we introduce the famous FKM isoparametric hypersurfaces in spheres and prove Theorem A. Sect. 3.1 is devoted to a review of Lawlor’s curvature criterion; for this part, we also refer the reader to \cite{TZ20}, \cite{Z25} for more details and remarks on Lawlor's criterion.  In Sect. 3.2 we compute the second fundamental forms and estimate the vanishing angles, while Sect. 3.3 is dedicated to computing the normal radius. The proof of Theorem B is then presented in Sect. 3.4. In Sect. 4, we study minimal product cones. Finally in Sect. 5, we give some comments and questions on FKM isoparametric foliations in $\mathbb{S}^n \times \mathbb{S}^n$ and on their area-minimizing cones. 

\medskip\noindent

\textbf{Notations.} In this paper, we always consider the following isometric embeddings:
$$(M, \nabla) \hookrightarrow \left(\mathbb{S}^n(1) \times \mathbb{S}^n(1), \bar{\nabla}\right) \hookrightarrow \left(\mathbb{S}^{2n+1}(\sqrt{2}), \widetilde{\nabla}\right) \hookrightarrow \left(\mathbb{R}^{2n+2}, D\right)$$
where $\nabla, \bar{\nabla}, \widetilde{\nabla}$, and $D$ denote the Levi-Civita connections associated with the induced metrics, respectively. $M$ is an orientable hypersurface of $\mathbb{S}^n \times \mathbb{S}^n$ with a unit normal vector field $N$. We also use $\nabla, \bar{\nabla}, \widetilde{\nabla}$, and $D$ to denote the associated gradient operators for functions, and let $\Delta, \bar{\Delta}, \widetilde{\Delta}$,  and $\Delta^{E}$ denote the associated Laplacian operators, respectively.

\medskip\noindent
\medskip\noindent

\section{Isoparametric hypersurfaces in $\mathbb{S}^n \times \mathbb{S}^n$}
\subsection{Level set hypersurfaces in $\mathbb{S}^n \times \mathbb{S}^n$}

\medskip\noindent

Under the natural isomorphism $T_{(x,y)}\left(\mathbb{S}^n \times \mathbb{S}^n\right) \cong  T_{x}\mathbb{S}^n \times T_{y}\mathbb{S}^n$, there is a product structure $P: T\left(\mathbb{S}^n \times \mathbb{S}^n\right) \rightarrow T\left(\mathbb{S}^n \times \mathbb{S}^n\right)$ such that
$$
P\left(x_1, x_2\right)=\left(x_1,-x_2\right), \quad \forall x_1, x_2 \in T \mathbb{S}^n,
$$
it is also referred to as almost product structure in some papers, as in \cite{LO74} and \cite{Chen21}. The tangent bundle $T\left(\mathbb{S}^n \times \mathbb{S}^n\right) \subset\left(\mathbb{S}^n \times \mathbb{S}^n\right) \times \mathbb{R}^{2n+2}$ has a natural decomposition with respect to the decomposition $\mathbb{R}^{2n+2}=\mathbb{R}^{n+1} \times \mathbb{R}^{n+1}$ and the natural projections $\pi_i: \mathbb{R}^{2n+2} \rightarrow \mathbb{R}^{n+1}, \text { with } \pi_i(x_1,x_2)=x_i \text { for } i=1,2$. The standard theory for Riemannian connections of products of Riemannian manifolds (cf. \cite[p.139]{Car92}) implies that $\bar{\nabla} P=0$. Moreover, $P$ preserves the metric. The geometries of hypersurfaces in $\mathbb{S}^{n} \times \mathbb{S}^{n}$ was studied in \cite{LO74}, \cite{MT13}, \cite{Ur19}, \cite{Chen21}, \cite{GHMY22}, \cite{LWW23}, \cite{LVWY24}, etc, see the book \cite{Ur25} for more information.

\medskip\noindent

As a product manifold, the Riemannian curvature tensor ${\rm R}$ of $\mathbb{S}^n \times \mathbb{S}^n$ is:
$$
\begin{aligned}
& {\rm R}\left(\left(x_1, x_2\right),\left(y_1, y_2\right),\left(z_1, z_{2}\right),\left(w_1, w_2\right)\right) \\
= & {\rm R}\left(x_1, y_1, z_1, w_1\right)+{\rm R}\left(x_2, y_2, z_2, w_2\right) \\
= & \langle x_1, z_1\rangle \langle y_1, w_1\rangle-\langle x_1, w_1\rangle \langle y_1, z_1\rangle 
 +\langle x_2, z_2\rangle \langle y_2, w_2\rangle-\langle x_2, w_2\rangle \langle y_2, z_2\rangle
\end{aligned}
$$
where $X=\left(x_1, x_2\right) \in \Gamma \left( T\left(\mathbb{S}^n \times \mathbb{S}^n\right)\right) $ (resp. $Y=\left(y_1, y_2\right)$, $Z=\left(z_1, z_2\right)$ and $W=\left(w_1, w_2\right)$), with respect to the product structure $P$, it can be also written as
$$
\begin{aligned}
{\rm R}(X, Y, Z, W)=\frac{1}{2}\{& \langle X, Z\rangle \langle Y, W \rangle+\langle PX, Z\rangle \langle PY, W \rangle\\
& -\langle X, W\rangle \langle Y, Z \rangle-\langle PX, W\rangle \langle PY, Z \rangle  \}.
\end{aligned}
$$

\medskip\noindent

Direct computations show that $\mathbb{S}^n \times \mathbb{S}^n$ is Einstein manifold and it has constant Ricci curvature $\operatorname{Ric}_{\mathbb{S}^n \times \mathbb{S}^n}=n-1$ and constant scalar curvature ${\rm S}_{\mathbb{S}^n \times \mathbb{S}^n}=2n(n-1)$. 

\medskip\noindent

With respect to the product structure $P$, F. Urbano \cite{Ur19} defined the following product angle function
$$
C:= \langle P N, N \rangle
$$
which characterizes the differences for the two projections of $N$. It is known that a hypersurface $M$ of $\mathbb{S}^{n} \times \mathbb{S}^{n}$ satisfying $C^2 \equiv1$ is locally an open part of $M' \times \mathbb{S}^{n}$ or $\mathbb{S}^{n} \times M''$ (cf. \cite[Lemma 2.3]{Chen21} and references therein), where $M', M''$ are hypersurfaces in $\mathbb{S}^{n}$. Moreover, one can set (cf. \cite{LWW23}) 
$$
\cos s=\left|\pi_1(N)\right|, \sin s=\left|\pi_2(N)\right|
$$
for some continuous function $s: M \rightarrow[0, \pi / 2]$, then $C=\cos 2 s$, $C \equiv \pm 1$ is equivalent to $s \equiv 0$ or $\pi / 2$, and $C \equiv 0$ is equivalent to $\left|\pi_1(N)\right|=\left|\pi_2(N)\right|$. We note here that $s$ and $C$ may be non-constant, see the examples in \cite[Sect. 3]{LWW23}.

\medskip\noindent

We now study the level set hypersurface $M$ of $\mathbb{S}^n \times \mathbb{S}^n$ given by the regular level set of some smooth function $F$ on $\mathbb{R}^{n+1} \times \mathbb{R}^{n+1}$. The unit normal vector $N(x,y)$ at $(x,y) \in M$ is
$$
N(x,y)=\frac{\bar{\nabla} F}{|\bar{\nabla} F|}(x,y)=\frac{DF- \langle DF,x \rangle x-\langle DF,y \rangle y}{|DF- \langle DF,x \rangle x-\langle DF,y \rangle y|},
$$
where we identify $x$ with $(x,0)$, $y$ with $(0,y)$ respectively.

\medskip\noindent

Let $N_{1}:=\left(\frac{ \partial F}{\partial x}-\langle DF,x \rangle x\right)$ and $N_{2}:=\left(\frac{ \partial F}{\partial y}-\langle DF,y \rangle y\right)$, then $\pi_{1}N=\frac{N_{1}}{|\bar{\nabla} F|}$ and $ \pi_{2}N=\frac{N_{2}}{|\bar{\nabla} F|}$, the Jordan angle function $\cos \theta=\frac{|N_{1}|}{|\bar{\nabla} F|}$ and $\sin \theta=\frac{|N_{2}|}{|\bar{\nabla} F|}$. Then Urbano's $C$ function 
$$
C=\langle PN,N \rangle=\cos 2\theta =\frac{|N_{1}|^2-|N_{2}|^2}{|\bar{\nabla} F|^2}=\frac{|N_{1}|^2-|N_{2}|^2}{|N_{1}|^2+|N_{2}|^2}.
$$

\medskip\noindent

In case $C^2<1$, $ N_{1},  N_2$  are both non-zero vectors, and the following canonical frame on $\mathbb{S}^n \times \mathbb{S}^n$ can be well-defined along  $M$:
\begin{equation}\label{F}
    \begin{aligned}
      e_{0}=\frac{1}{|N_{1}|}(&N_{1},0),e_{1}=(\varepsilon_{1},0),\ldots,e_{n-1}=( \varepsilon_{n-1},0),\\ e_{0}^{\prime}=\frac{1}{|N_{2}|}(& 0,N_{2}),e_{1}^{\prime}=(0,\varepsilon_{1}^{\prime}), \ldots,e_{n-1}^{\prime}=(0,\varepsilon_{n-1}^{\prime}),\\ 
    \end{aligned}
\end{equation}
where $\{x,\varepsilon_{0}:=\frac{N_{1}}{|N_{1}|}, \varepsilon_{1},\ldots,\varepsilon_{n-1}\}\text{ and }\{y,\varepsilon_{0}^{\prime}:=\frac{N_{2}}{|N_{2}|},\varepsilon_{1}^{\prime},\ldots,\varepsilon_{n-1}^{ \prime}\} $
are orthonormal frames of $\mathbb{R}^{n+1}$ defined at $x$ and $y$ respectively, and $N=\cos \theta e_{0}+ \sin \theta e_{0}^{\prime}.$

\medskip\noindent 

If $C^2<1$, an important observation is that, for $(x,y)\in M$, the two slices of $M$, $M_{1}$: $M \cap\left(\mathbb{S}^{n} \times\{y\}\right) \hookrightarrow \mathbb{S}^{n} \times\{y\}$  and $M_{2}$: $\left(\{x\}\times \mathbb{S}^{n}\right) \cap M \hookrightarrow\{x\} \times \mathbb{S}^{n}$ are both regular level sets of the further  restrictions of $F$ when $y$ (or $x$) is fixed, then $\varepsilon_{0}$ and $\varepsilon_{0}^{\prime}$ are just their unit normal vector fields at $(x,y)$ in the unit spheres under natural diffeomorphism: $\mathbb{S}^{n} \times\{y\} \cong \mathbb{S}^{n}$ and $\{x\} \times \mathbb{S}^{n} \cong \mathbb{S}^{n}$. In case $C^2<1$, there is an extra tangent vector of $M$:
$$
(PN)^{T}=PN-CN=\left(\frac{1-C}{|\bar{\nabla} F|}N_{1},-\frac{1+C}{|\bar{\nabla} F|}N_2 \right),
$$
with $|(PN)^{T}|=1-C^2$, then 
$$
T:=\frac{PN-CN}{|PN-CN|} =\left(\sqrt{\frac{1-C}{2}}e_{0},-\sqrt{\frac{1+C}{2}}e_{0}^{\prime} \right)
$$
is a local unit tangent vector field. 

\medskip\noindent

However, one can check directly that $T$ can be defined for any case. If $C(p)=1$, then $N(p)$ is located in the first entry, all orthogonal vectors to $y$ are also tangent vectors of $M$, by setting $\varepsilon_{0}^{\prime}$ to be one such orthogonal vector at $y$, $T$ is still a tangent vector. Similar for the case $C(p)=-1$.

\medskip\noindent

\textbf{Definition 2.1.} \textit{A canonical frame on level set hypersurfaces in $ \mathbb{S}^n \times \mathbb{S}^n$ is give as:
\begin{equation}\label{cf}
    \{\tau_1=e_{1},\cdots,\tau_{n-1}=e_{n-1},\tau_{n}=e_{1}^{\prime},\cdots,\tau_{2n-2}=e_{n-1}^{\prime},\tau_{2n-1}=T \}.
\end{equation}
The position vector is denoted by $X=(x,y)$. Particularly, if $C\equiv0$, $T=\frac{e_{0}-e_{0}^{\prime}}{\sqrt{2}}$ and the unit normal is $N=\frac{e_{0}+e_{0}^{\prime}}{\sqrt{2}}$.}

\medskip\noindent

The mean curvature of a level set hypersurface can be computed.

\medskip\noindent

\textbf{Lemma 2.2. (see \cite[Theorem 3.3]{CR15})} \textit{Assume $M$ is a regular level set of a smooth function $F$ on a Riemannian manifold, denote the gradient of $F$ by $\bar{\nabla} F$, then the mean curvature of $M$ is:
$$
H=\frac{\langle\bar{\nabla}|\bar{\nabla} F|^{2},\bar{\nabla}F\rangle}{2|\bar{\nabla} F|^{3}}-\frac{\bar{\Delta} F}{|\bar{\nabla} F|}.
$$}

\medskip\noindent

Moreover, we perform a computation for the length of the second fundamental forms and give the following formula, which is important in later applications.

\medskip\noindent

\textbf{Theorem 2.3.} \textit{Assume $M $ is a level set hypersurface given by a regular level set of some function $F$ on $\mathbb{S}^n \times \mathbb{S}^n$, then the square norm of the second fundamental forms of $M \subset \mathbb{S}^n \times \mathbb{S}^n$ is
\begin{equation}\label{2nd4}
    ||B||^2=\frac{1}{|\bar{\nabla} F|^2} \left( \frac{1}{2} \bar\Delta|\bar\nabla F|^2- \langle\bar\nabla F, \bar\nabla \bar\Delta F\rangle - (n-1)|\bar{\nabla} F|^2 -\frac{\left|\bar{\nabla}\left|\bar{\nabla} F\right|^2\right|^2}{2|\bar{\nabla} F|^2}+\frac{\left\langle \bar{\nabla}\left|\bar{\nabla} F\right|^2, \bar{\nabla} F \right\rangle^2 }{4|\bar{\nabla} F|^4}\right),
\end{equation}}.

\medskip\noindent

\textbf{Proof:}   For $M \subset \mathbb{S}^n \times \mathbb{S}^n$ a regular level set of an isoparametric function $F$ on $\mathbb{S}^n \times \mathbb{S}^n$, its second fundamental forms are: $B\left(f_i, f_j\right)=\left\langle\nabla_{f_i} f_j, \frac{\bar{\nabla} F}{|\bar{\nabla} F|}\right\rangle =-\frac{H F(f_{i}, f_{j})}{|\bar{\nabla} F|}$, where $\{f_{i}\}$ are orthonormal frame for the tangent bundle of $M$. Thus
\begin{equation}\label{2nd1}
\begin{aligned}
|\bar{\nabla} F|^2||B||^2& =\sum_{i j} \operatorname{Hess} F\left(f_i, f_j\right)^2=\left\langle \bar{\nabla}_{f_i} \bar{\nabla} F, f_j\right\rangle^2\\
&=\sum_i\left[\left\langle\bar{\nabla}_{f_i} \bar{\nabla} F, \bar{\nabla}_{f_i} \bar{\nabla} F \right\rangle^2-\left\langle\bar{\nabla}_{f_i} \bar{\nabla} F, \frac{\bar{\nabla} F}{|\bar{\nabla} F|}\right\rangle^2\right] \\
& =\left(||\operatorname{Hess} F||^2-\frac{\left|\bar{\nabla}_{\bar{\nabla} F} {\bar{\nabla} F}\right|^2}{|\bar{\nabla} F|^2}\right)-\sum_{i} \frac{1}{|\bar{\nabla} F|^2} \left\langle \frac{\bar{\nabla}\left|\bar{\nabla} F\right|^2}{2}, f_i\right\rangle^2 \\
&=||\operatorname{Hess} F||^2-\frac{\left|\bar{\nabla}\left|\bar{\nabla} F\right|^2\right|^2}{4|\bar{\nabla} F|^2}-\frac{1}{|\bar{\nabla} F|^2}\left( \left| \frac{\bar{\nabla}\left|\bar{\nabla} F\right|^2}{2}\right|^2-\left\langle \frac{\bar{\nabla}\left|\bar{\nabla} F\right|^2}{2}, \frac{\bar{\nabla} F}{|\bar{\nabla} F|} \right\rangle^2 \right)\\
&=||\operatorname{Hess} F||^2-\frac{1}{2|\bar{\nabla} F|^2}\left|\bar{\nabla}\left|\bar{\nabla} F\right|^2\right|^2+\frac{1}{4|\bar{\nabla} F|^4}\left\langle \bar{\nabla}\left|\bar{\nabla} F\right|^2, \bar{\nabla} F \right\rangle^2
\end{aligned}
\end{equation}
where $\bar{\nabla}_{\bar{\nabla} F} {\bar{\nabla} F}=\frac{1}{2} \bar{\nabla}\left|\bar{\nabla} F\right|^2$ can be obtained by the symmetric properties of $\operatorname{Hess} F$ (for example, see \cite[p.61]{CM11}), the above formula also holds in case $\mathbb{S}^n \times \mathbb{S}^n$ is replaced by a general Riemannian manifold.

\medskip\noindent

The Hessian term of $F$ can be computed by the Bochner formula (for example, see \cite[Prop. 1.47]{CM11}):
\begin{equation}\label{2nd2}
\frac{1}{2} \bar\Delta|\bar\nabla F|^2=\|\operatorname{Hess} F\|^2+\langle\bar\nabla F, \bar\nabla \bar\Delta F\rangle+\operatorname{Ric}_{\mathbb{S}^n \times \mathbb{S}^n}(\bar\nabla F, \bar\nabla F).
\end{equation}

\medskip\noindent

We now compute the term of Ricci curvature. By frame \eqref{cf}, and set $\bar{\nabla} F=|N_{1}|e_{0}+|N_{2}|e_{0}^{\prime}:=a e_0+b e_0^{\prime}$, then
\begin{equation}\label{2nd3}
\begin{aligned}
\operatorname{Ric}( \bar{\nabla} F,\bar{\nabla} F)=&\sum_A R\left(E_A, \bar{\nabla} F, E_A,\bar{\nabla} F\right)\\
=& R\left(e_0, a e_0+b e_0^{\prime}, e_0, a e_0+b e_0^{\prime}\right)+R\left(e_0^{\prime}, a e_0+b e_0^{\prime}, e_0^{\prime}, a e_0+b e_0^{\prime}\right) \\
&+R\left(e_i, a e_0+b e_0^{\prime}, e_i, a e_0+b e_0^{\prime}\right)+R\left(e_i^{\prime}, a e_0+b e_0^{\prime}, e_i^{\prime}, a e_0+b e_0^{\prime}\right) \\
= & a^2 \sum_{i=1}^{n-1} R\left(e_i, e_0, e_i, e_0\right)+b^2 \sum_{i=1}^{n-1} R\left(e_i^{\prime}, e_0^{\prime}, e_i^{\prime}, e_0^{\prime}\right) \\
= & \sum_{i=1}^{n-1} a^2+\sum_{i=1}^{n-1} b^2 \\
=&(n-1)|\bar{\nabla} F|^2.
\end{aligned} 
\end{equation}

\medskip\noindent

Combining \eqref{2nd1},\eqref{2nd2} and \eqref{2nd3}, we finally get the formula \eqref{2nd4}.
$\Box$

\medskip\noindent

\textbf{Corollary 2.4.} \textit{An isoparametric hypersurface in $\mathbb{S}^n \times \mathbb{S}^n$ has constant $||B||^2$.}

\medskip\noindent

\subsection{FKM isoparametric hypersurfaces in $\mathbb{S}^n \times \mathbb{S}^n$ and the proof of Theorem A}

\medskip\noindent

Given an Riemannian manifold $X$, a non-constant smooth function $f$ on $X$ is called isoparametric if there exist smooth functions $a,b$ such that
\begin{equation}\label{oc1}
  |\nabla f|^2=b(f) \ {\rm and} \ \Delta f=a(f),  
\end{equation}
where $\nabla$ and $\Delta$ are the gradient operator and Laplacian operator on $X$, respectively. A regular level set $f^{-1}(t)$ is called an isoparametric hypersurface of $X$. For space forms, a hypersurface $M$ is isoparametric if and only if it has constant principal curvatures. Isoparametric hypersurfaces in $X=\mathbb{R}^l$ or $X=\mathbb{H}^l$ were classified by Cartan and Segre in a very early time. For $X=\mathbb{S}^l$, let $g$ denote the number of distinct principal curvatures of $M$, and let $m_{\alpha}(1 \leq \alpha \leq g)$ denote their multiplicities. The structure theorems of Münzner \cite{M80} stated that $g$ can only be $1,2,3,4,6$ and $m_{\alpha}=m_{\alpha+2} \ ({\rm mod} \ g)$. Classifying isoparametric hypersurfaces in spheres is one of the most difficult and fascinating problems in differential geometry.

\medskip\noindent

The final case—and also the most challenging one—in the classification of isoparametric hypersurfaces in the sphere is when $g=4$. Thanks to recent elegant results by Cecil-Chi-Jensen \cite{CCJ07} and Chi \cite{Chi20}, the classification is now finally complete. Such isoparametric hypersurfaces must be either OT-FKM type or belong to two known homogeneous examples; for further details, see the book \cite[Sect. 3.9]{CR15} and the excellent survey \cite{Chi18}, \cite{GQTY25}. Here, the "OT-FKM type" means that the level sets
$$
\textbf{M}_t=\mathbf{F}^{-1}(t) \bigcap S^{l}(1), \ \ \ -1\leq t \leq 1
$$
where
$$
\mathbf{F}(\widetilde{x})= \langle \widetilde{x}, \widetilde{x} \rangle^2-2 \sum_{i=0}^{m} \langle P_i \widetilde{x}, \widetilde{x} \rangle^2, \ \ \ m\geq 1
$$
form the leaves of the isoparametric foliation in the unit sphere. $\mathbf{F}$ is the Cartan-Münzner polynomial for the examples constructed by Ferus, Karcher, and Münzner (\cite{FKM81}), which generalize the examples of Ozeki-Takeuchi (\cite{OT75}, \cite{OT76}). We briefly refer to $\mathbf{F}$ as \textit{FKM isoparametric polynomials}. $P_i (0\leq i \leq m)$ are a family of symmetric orthogonal transformations, i.e.,
$$
\left\{\begin{array}{c}
P_i P_j+P_j P_i=2 \delta_{ij} I \\
P_i=P_{i}^{T}
\end{array}\right.
$$
which are referred to as a \textit{symmetric Clifford system} or \textit{Clifford system} (cf. \cite[Sect. 3.9]{CR15} and \cite[Sect. 5.4]{BW03}). Such a structure exists if and only if $l+1=2(n+1)$ for some $n$ and $n+1=k\delta(m)$, $\delta(m)$ denotes the dimension of an irreducible module of Clifford algebra $Cl_{m-1}$. Their multiplicities satisfy $(m_{1},m_{2})=(m,n-m)$. All possible values of $\delta(m)$ are the following (see, for example, \cite[Sect. 5.4]{BW03}, \cite[p.525]{TY13}, \cite[Sect. 3.9]{CR15})

\begin{equation}\label{table}
\begin{array}{|c|c|c|c|c|c|c|c|c|c|c|}
\hline m & 1 & 2 & 3 & 4 & 5 & 6 & 7 & 8 & \cdots & m+8 \\
\hline \delta(m) & 1 & 2 & 4 & 4 & 8 & 8 & 8 & 8 & \cdots &16 \delta(m) \\
\hline
\end{array}
\end{equation}
and by \cite[Theorem 3.72]{CR15} 
\begin{equation}\label{bmi}
   \left\{
\begin{aligned}
    |\widetilde{\nabla} \mathbf{F}|^2&=16(1-F^2)\\
    \widetilde{\Delta}\mathbf{F}&=8(n-2m)-8(n+2)F
\end{aligned}
   \right.
\end{equation}
where $\widetilde{\nabla}, \widetilde{\Delta}$ denote the Riemannian connection and the Laplacian of $\mathbb{S}^{2n+1}(1)$.

\medskip\noindent

Isoparametric hypersurfaces in general Riemannian manifolds were studied in \cite{W87}, \cite{TT95}, \cite{GT09}, \cite{Mi12}, \cite{GT14}, \cite{GTY15}, \cite{QT16}, etc. Recently, people concern on isoparametric
hypersurfaces in the products of two real space forms (for the newest progress, see \cite{DP24}, \cite{TXY25}, \cite{DP25}). For the ambient manifold $\mathbb{S}^2 \times \mathbb{S}^2$, Urbano \cite{Ur19} proved that the product angle function $C$ is a constant function and finally classified the isoparametric hypersurfaces in $\mathbb{S}^2 \times \mathbb{S}^2$ through a classification of hypersurfaces with constant $C$(\cite[Corollary 3]{Ur19}): they must be either (i) $\mathbb{S}^1(r) \times \mathbb{S}^2$, $r\in (0,1]$; or (ii) $M_{t}=\{(p,q)\in \mathbb{S}^2 \times \mathbb{S}^2| \langle p,q \rangle=t\}$, $t\in (-1,1)$. The former case is trivial, while the latter case is a special instance of \textit{the splitting form of FKM isoparametric polynomials} which can be explicitly described below.

\medskip\noindent

For the above symmetric Clifford system, up to an orthogonal transformation, one can always express it as follows (cf. \cite[Sect. 3.9]{CR15})
$$
P_0=\left(\begin{array}{rr}
I_l & 0 \\
0 & -I_l
\end{array}\right), \quad P_1=\left(\begin{array}{cc}
0 & I_l \\
I_l & 0
\end{array}\right), \quad P_{1+\alpha}=\left(\begin{array}{rr}
0 & A_\alpha \\
-A_\alpha & 0
\end{array}\right) \quad \text { for } 1 \leq \alpha \leq m-1,
$$
where $\left\{A_1, A_2, \ldots, A_{m-1}\right\}$ are skew-symmetric representations of the Clifford algebra $Cl_{m-1}$ on $\mathbb{R}^{n+1}$, i.e.
\begin{equation}\label{skc}
    \left\{\begin{array}{c}
A_p A_q+A_q A_p=-2 \delta_{p q} I \\
A_p=-A_p^{T}.
\end{array}\right.
\end{equation}
Up to a scaling, we can set $\widetilde{x}= (x,y) \in \mathbb{S}^{2n+1}(\sqrt{2})$, where $x,y \in \mathbb{S}^n(1)$, and it can be checked that (as did in \cite[Equ. 4.13]{GT23} under a different research background)
\begin{equation}\label{redu}
    \mathbf{F}(\widetilde{x})=(|x|^2+|y|^2)^2-2(|x|^2-|y|^2)^2-8\langle x, y\rangle^2-8\sum_{\alpha=1}^{m-1}\langle A_\alpha x, y\rangle^2,
\end{equation}
the restriction of $\mathbf{F}$ on $\mathbb{S}^{n} \times \mathbb{S}^{n} \subset \mathbb{S}^{2n+1} $ is equivalent to the restriction of
\begin{equation}
    F(x,y)=\langle x, y\rangle^2+\sum_{\alpha=1}^{m-1}\langle A_\alpha x, y\rangle^2
\end{equation}
on $\mathbb{S}^{n} \times \mathbb{S}^{n}$, and $F$ is a bi-homogeneous isoparametric polynomial with respect to $x,y$. 

\medskip\noindent

\textbf{Definition 2.5.} \textit{We call $F(x,y)=\langle x, y\rangle^2+\sum_{\alpha=1}^{m-1}\langle A_\alpha x, y\rangle^2$ the splitting form of FKM isoparametric polynomials on $\mathbb{R}^{n+1} \times \mathbb{R}^{n+1}$.}

\medskip\noindent

In the special case $m=1$, its splitting form yields the isoparametric polynomial $\langle x, y\rangle$ and associated isoparametric hypersurfaces $M_{t}$ in  $\mathbb{S}^{n} \times \mathbb{S}^{n}$ as given by \cite[Theorem 1]{QT16} and F. Urbano \cite[Remark 1]{Ur19}. 

\medskip\noindent

We now make it precise that

\medskip\noindent

\textbf{Theorem A} \textit{The restrictions of $F=\langle x, y\rangle^2+\sum_{q=1}^{m-1}\langle A_q x, y\rangle^2$  to $\mathbb{S}^{n} \times \mathbb{S}^{n}$ give rise to isoparametric hypersurfaces which satisfying $C\equiv 0$.}

\medskip\noindent

\textbf{Proof:}  (1) $ \frac{1}{2}F_{x_i}=\langle x, y\rangle y_i+\sum_{q}\left\langle A_q x, y\right\rangle y^{T} A_q e_i$, then by applying \eqref{skc}
\begin{equation}\label{czero}
\begin{aligned}
\frac{1}{4}\sum_{i} F_{x_i}^2 =& \langle x, y\rangle|y|^2+2\langle x, y\rangle \sum_q \langle A_q x, y \rangle\left(y^{T} A_q y\right) \\
&+\sum_{p,q}\langle A_p x, y\rangle \langle A_q x, y\rangle \sum_i\left(y^{T} A_p e_i y^{T} A_q e_i\right) \\
&=\langle x, y\rangle|y|^2+\sum_{q}\langle A_q x, y\rangle^2|y|^2 \\
&=F|y|^2.
\end{aligned}
\end{equation}
Similarly for $\sum_{i} F_{y_i}^2$, then one gets that $|DF|^2=4F(|x|^2+|y|^2)$. Since the normal component of $D F$ is $2 F \cdot (x, 0)+2 F \cdot (0, y)$, then we find that $|\bar{\nabla} F|^2=8F-8F^2$.

\medskip\noindent

(2) $ \frac{1}{2} F_{x_i x_i}=y_i^2+\sum_{q}\left(y^{T} A_q e_i\right)^2$ implies that
\begin{equation}
\frac{1}{2} \sum_i F_{x_i x_i}=|y|^2+\sum_q\left| y \right|^2= m | y|^2
\end{equation}
similarly for $\sum_{i} F_{y_i y_i}$, then one gets that $\Delta^{E} F=2m(|x|^2+|y|^2)$. To compute $\bar \Delta F$, we choose orthonormal frame  $$\left\{e_{0}=(x,0), e_{1}=(\varepsilon_{1},0), \ldots, e_{n}=(\varepsilon_{n},0),e_{0}^{\prime}=(0,y), e_{1}=(0,\varepsilon_{1}^{\prime}), \ldots, e_{n}=(0,\varepsilon_{n}^{\prime})\right\}$$ on some point $(x, y) \in \mathbb{S}^n \times \mathbb{S}^n$, such that  $\varepsilon_{1}, \ldots, \varepsilon_{n}$  are tangent to  $\mathbb{S}^n$ at $x$, $\varepsilon_{1}^{\prime}, \ldots, \varepsilon_{n}^{\prime}$  are tangent to  $\mathbb{S}^n$ at $y$. From the standard theory for Riemannian connections of products of manifolds, we have $e_{0}(F)=\langle DF,e_{0}\rangle=2F$, $D_{e_{0}}e_{0}=e_{0}$ and $\langle D_{e_{i}}e_{i} , e_{0}\rangle=-\langle D_{e_{i}}e_{0} , e_{i}\rangle=-\langle e_{i}, e_{i}\rangle =-1$, thus
\begin{equation}
    \begin{aligned}
    \Delta^{E} F&=
       e_0e_0(F)+e_0^{\prime}e_0^{\prime}(F)+e_ie_i(F)+e_i^{\prime}e_i^{\prime}(F)-D_{e_{0}}e_{0}(F)-D_{e_{0}^{\prime}}e_{0}^{\prime}(F)-D_{e_{i}}e_{i}(F)-D_{e_{i}^{\prime}}e_{i}^{\prime}(F) \\
       &=\bar{\Delta} F+4F+4F-2F-2F+2nF+2nF \\
       &=\bar{\Delta} F+4(n+1) F,
    \end{aligned}
\end{equation}
then $\bar{\Delta}F=4m-4(n+1)F.$

\medskip\noindent

In summary, we have
\begin{equation}
   \left\{
\begin{aligned}
    |\bar{\nabla} F|^2&=8F-8F^2\\
    \bar{\Delta}F&=4m-4(n+1)F
\end{aligned}
   \right.
\end{equation}
then, $F|_{\mathbb{S}^n \times \mathbb{S}^n}$ is an isoparametric function which takes values $[0,1]$. This whole isoparametric foliation is the restriction of the FKM isoparametric foliation in $\mathbb{S}^{2n+1}(\sqrt{2})$ to $\mathbb{S}^n \times \mathbb{S}^n$.

\medskip\noindent

To prove Urbano's $C$ function for the level sets: $F^{-1}(t)(0<t<1)$ are $C=0$, note that $|N_{1}|^2=\left|\frac{ \partial F}{\partial x} \right|^2-4F^2$ and $|N_{2}|^2=\left|\frac{ \partial F}{\partial y} \right|^2-4F^2$, it is equivalent to prove that $\left|\frac{ \partial F}{\partial x} \right|^2=\left|\frac{ \partial F}{\partial y} \right|^2$. Following \eqref{czero}, restricting on $\mathbb{S}^n \times \mathbb{S}^n$, $\left|\frac{ \partial F}{\partial x} \right|^2=\left|\frac{ \partial F}{\partial y} \right|^2=4 F $.
 $\Box$
 
\medskip\noindent

\textbf{Remark 2.6.} \textit{To maintain consistency with conventional notation, we may renormalize $F$ to $\widetilde{F}=|x|^2|y|^2-2F$ such that the focal variety associated to the global maximum value of $\widetilde{F}$ is then $M_{+}=\widetilde{F}^{-1}(1)=F^{-1}(0)$, $M_{-}=\widetilde{F}^{-1}(-1)=F^{-1}(1)$.  As in \cite[p. 174]{CR15}, the focal variety $M_{+}=F^{-1}(0)$ is still the focal variety $\mathbf{F}^{-1}(4)$ (under a rescaling of $\sqrt{2}$) of the original FKM isoparametric foliations, they are also called the Clifford-Stiefel manifolds. However, the focal variety $M_{-}=F^{-1}(1)$ is the intersection of $\mathbf{F}^{-1}(-4)$ with $\mathbb{S}^n(1) \times \mathbb{S}^n(1)$.}

\medskip\noindent

Following Lemma 2.2, direct computation shows that $\Sigma:=F^{-1}(\frac{m-1}{n-1})$ are minimal isoparametric hypersurfaces in $\mathbb{S}^n \times \mathbb{S}^n$. Then by Lemma 1.2, $\Sigma$ is also a minimal submanifold in $\mathbb{S}^{2n+1}(\sqrt{2})$. 

\medskip\noindent

\textbf{Remark 2.7.}  \textit{By Lemma 2.2 and \eqref{bmi}, the minimal FKM isoparametric hypercones in $\mathbb{R}^{2n+2}$ is given by $\mathbf{F}(\widetilde{x})=\frac{n-2m}{n}|\widetilde{x}|^4$. Following \eqref{redu}, its intersection with Simons cones are 
$$
\langle x, y\rangle^2+\sum_{q=1}^{m-1}\langle A_q x, y\rangle^2=\frac{m}{n}|x|^2|y|^2 \ \ \ {\rm and} \ \ \ |x|=|y|,
$$
it is not $\Sigma$. Since by counting the multiplicity as in \cite[Theorem 3.72]{CR15}, it must $n-m>0$.}

\medskip\noindent

In the process of the computations, one can directly verify that

\medskip\noindent

\textbf{Proposition 2.8.} \textit{The restriction of $F=\langle x, y\rangle^2+\sum_{q=1}^{m-1}\langle A_q x, y\rangle^2$ on general $\mathbb{S}^n(a) \times \mathbb{S}^n(b)$ also gives isoparametric hypersurfaces.}

\medskip\noindent

The FKM isoparametric hypersurfaces in general $\mathbb{S}^n(a) \times \mathbb{S}^n(b)$ have another explanation. In \cite{W86}, Wang gives a definition for isoparametric maps (also see \cite{W94} and \cite{X93}), which is more general than the well-known definition used by Terng \cite{T85}. Wang's definition is a direct generalization of \eqref{oc1}, and it ignores the integrability condition of the horizontal distribution. This definition is useful in studying minimal submanifolds or harmonic maps, see \cite[Chapter 4]{ER93}.

\medskip\noindent

\textbf{Definition 2.9 (\cite{W86}).} \textit{Let $X$ be a smooth complete Riemannian manifold. A smooth map $f: X \rightarrow \mathbb{R}^k$ is called isoparametric if
$$
\begin{aligned}
\left\langle d f_i, d f_j\right\rangle & =a_{i j}\left(f_1, \ldots, f_k\right) \\
\Delta f_i & =b_i\left(f_1, \ldots, f_k\right)
\end{aligned}
$$
for $i, j=1,2, \ldots, k$, where $a_{i j}$ are smooth function on $f(X) \subset \mathbb{R}^k$.}

\medskip\noindent

Then the isoparametric hypersurfaces in $\mathbb{S}^n(a) \times \mathbb{S}^n(b)$ can be seen as the regular fibers of the following isoparametric map $\widetilde{F}: \mathbb{R}^{n+1} \times \mathbb{R}^{n+1} \rightarrow \mathbb{R}^{3}$ given by 
$$
f_1=F(x,y), f_2=|x|^2-|y|^2, f_3=|x|^2+|y|^2
$$
or the isoparametric map $\pi: \mathbb{S}^{2n+1} \rightarrow \mathbb{R}^{2}$ given by 
$$
f_1=F(x,y), f_2=|x|^2-|y|^2,
$$
it could be related to the study of singular Riemannian foliations in \cite{R14}, see the comments in Sect. 5 and an introduction in \cite{GQTY25}.

\medskip\noindent

\textbf{Remark 2.10.} \textit{Particularly, we have the following isoparametric examples
\begin{equation}
    F(X,Y)=|\langle X, Y \rangle_{\mathbb{K}}|^2, \ \ \ \ \ X,Y \in \mathbb{K}^{n+1}, \ \ \mathbb{K}=\mathbb{R},\mathbb{C},\mathbb{H}, \mathbb{O},
\end{equation}
where $
    \langle X, Y \rangle_{\mathbb{K}}= \sum_{i=1}^{n+1} X_i \cdot \bar{Y_{i}}$, and “$\cdot$” denote the multiplication in normed division algebras. We note here that all these expressions also appeared in \cite{GT23}, especially in \cite[Equ. 4.13]{GT23} under a different geometric background.}     

\medskip\noindent
\medskip\noindent

\section{Area-minimizing subcones of Simons cones}

\medskip\noindent 

Let $M^{k-1} \subset \mathbb{S}^l \subset \mathbb{R}^{l+1}$ be an $(k-1)$-dimensional oriented closed embedded submanifold (or rectifiable current without boundary) in the unit sphere. The cone over $M$ is
$$
\mathcal{C}=\{t x: t \in[0, \infty) \text { and } x \in M\},
$$
then the truncated cone $\mathcal{C}_{1}:=\mathcal{C} \cap \textbf{B}^{l+1}(1)$ has the boundary $\partial \mathcal{C}_{1}=M$. $M$ is called the \textit{link} of $\mathcal{C}$. A cone $\mathcal{C}$ is said to be area-minimizing if $\mathcal{C}_{1}$ has the least area among all integral currents with boundary $M$. Given a smooth boundary $M$, 
the cone $\mathcal{C}$ owns an isolated singularity at the origin, and such cones are called regular cones (cf. \cite{HS85}). By Lawlor's deep research \cite{La91}, the area-minimizations of regular minimal cones depend heavily on the differential geometry properties of their link $M$.

\medskip\noindent

\subsection{Lawlor's curvature criterion}

\medskip\noindent 

One can  refer to \cite{M02}, \cite{Z18}, \cite{TZ20}, \cite{CX25}, \cite{CC26}, \cite{Z25} for more details and applications on Lawlor's curvature criterion. For readers’ convenience, we still provide a review here. Lawlor's main results are stated in Theorem 3.5.

\medskip\noindent 

Let $\mathcal{C}=C(\Sigma) \subset \mathbb{R}^{l+1}$ be a $k$-dimensional regular minimal cone with isolated singularity at the origin, where $\Sigma $ is a smooth, minimal submanifold in $\mathbb{S}^{l}(1)$ and $k \geq 3$. Firstly, we recall Lawlor's definitions for \textit{normal wedge} and \textit{normal radius}:

\medskip\noindent 

\textbf{Definition 3.1. (cf. \cite[Definition 1.1.1,1.1.2,1.1.3]{La91})} \textit{\emph{(1)} For $p\in \Sigma$ and $\alpha \geq 0$, let $U_{p}(\alpha)$ be the union of all open normal geodesics (in the sphere) of length $\alpha$ from $p$. The \textit{normal wedge} $W_{p}(\alpha)$ is defined as the conical set over $U_{p}(\alpha)$.}

\textit{\emph{(2)} Let $N_{p}:= {\rm max} \{\alpha | W_{p}(\alpha)\cap \mathcal{C}= \overrightarrow{op}\}$ , where $\overrightarrow{op} \in  \mathcal{C}$ is the ray through $p\in \Sigma$, we call $N_{p}$ the \textit{normal radius} at $p$. Equivalently, $N_{p}$ is the shortest normal geodesic in the unit sphere starting from $p$ that intersects another point of $\Sigma$.}

\medskip\noindent

Choose a unit normal vector $v$ at $p$, denote the $\mathbf{2}$-dimensional subspace spanned by $\{\overrightarrow{op}, v\}$ of $\mathbb{R}^{n}$ by $\mathcal{L}_{p,v}$. On $\mathcal{L}_{p,v}$, Lawlor select a curve in polar coordinate $(r(\theta),\theta)$: $S_{p}=r(\theta) ( \cos  \theta \ \overrightarrow{op}+ \sin  \theta \ v)$, then extend $S_{p}$ homothetically to a map $\Phi$ from the cone $\mathcal{C}$ to the ambient Euclidean space. Lawlor want to find some function $r(\theta)$ such that $S_{p}$ can extend to infinity when $\theta$ approach to some fixed angle $\theta_{0}$ and also the $k$-dimensional Jacobian: $J_{k}\Phi \geq 1$, it imply that the resulting map $\Phi$ is area non-decreasing and $\Phi$ takes values at a normal wedge of $\mathcal{C}$ at $p$. And, if pointwisely the normal wedges $ W_{p}(\theta_{0}(p))$ are disjoint for any $p,q \in \Sigma, p\neq q$:
\begin{equation}\label{Normal radius}
 W_{p}(\theta_{0}(p))\cap W_{q}(\theta_{0}(q))=\varnothing, 
\end{equation}
the inverse map $\Pi:=\Phi^{-1}$ is well-defined in an angle neighborhood of the cone $W:= \bigcup_{p\in \Sigma} W_{p}(\theta_{0}(p))$, then define that $\Pi$ maps every point outside $W$ to the origin, it follows that $\Pi$ is area non-increasing. 

\medskip\noindent

\textbf{Definition 3.2. (\cite[Sect. 1.1]{La91})} \textit{We call the above area non-increasing retraction map $\Pi$: the \textit{Lawlor retraction map}.}

\medskip\noindent

The existence of the Lawlor retraction map (\cite[p.10-12]{La91}) is reduced to the following ordinary differential equation which we just call it \textit{Lawlor ODE}.

\medskip\noindent

\textbf{Theorem 3.3.} \textit{If
\begin{equation}\label{va1}
    \frac{dr}{d\theta}\leq r \sqrt{r^{2k}({\rm cos}\theta)^{2k-2}\left(\mathrm{inf}_{v}{\rm det}(I-{\rm tan} \theta A_v)\right)^{2}-1}
\end{equation}
for any normal direction $v$, then $J_{k}\Pi \leq 1$, where $A_v$ is the shape operator of $\Sigma \subset \mathbb{S}^{l}$ in normal direction $v$.}

\medskip\noindent

The inequality in \eqref{va1} can be replaced by the equality to make the possibility for the non-intersection of $\theta_0$-normal wedges bigger. Then, the solution is analyzed by Lawlor:  either $\frac{dr}{d\theta}$ vanish at some positive $\theta(p)$ or $r\rightarrow \infty$ as $\theta \rightarrow \theta_{0}(p)$, in the latter case, Lawlor can construct his retraction map $\Pi$ if in addition \eqref{Normal radius} is satisfied for the normal wedges with the angle $\theta_{0}$, finally it follows that the cone is area-minimizing (see \cite[Theorem 1.2.1]{La91}).

\medskip\noindent
 
\textbf{Definition 3.4. (\cite[Definition 1.1.7]{La91})} \textit{We call the smallest $\theta_{0}(p)$ among all the normal directions $v$ the vanishing angle at $p$.}

\medskip\noindent

The global condition \eqref{Normal radius}  depends on the submanifold geometry of $\Sigma$ and is not easy to check directly (see Remark 3.10). Lawlor gave the following simplified criterion related to the normal radius:

\medskip\noindent

\textbf{Theorem 3.5. (\cite[Theorem 1.3.5]{La91})} \textit{If the vanishing angle $\theta_{0}(p)$ exists for any points $p\in \Sigma$ and pointwisely satisfy
\begin{equation}\label{two times vanishing angle}
2\theta_{0}(p)\leq N(p),
\end{equation}
where $N(p)$ is the normal radius at $p$, then the cone $\mathcal{C}=C(\Sigma)$ is area-minimizing (in the sense of mod $2$ when $\Sigma$ is nonorientable).}

\medskip\noindent

By a vanishing calibration argument, Lawlor showed that the inequality \eqref{va1} is equivalent to the following:
\begin{equation}\label{new ODE}
   \left(g(t)-\frac{t}{k} g^{\prime}(t)\right)^2+\left(\frac{g^{\prime}(t)}{k}\right)^2 \leq p^2(t),
\end{equation}
where $p(t):=\mathrm{inf}_{v}{\rm det}(I-{\rm tan} \theta A_v)$. The associated transformation is $g(t)=(r \cos \theta)^{-k}, t=\tan \theta$. Thus, if the inequality \eqref{new ODE} holds for some solutions $g(t)$
hitting $0$ soon enough, i.e., satisfying $g(0)=1$, and $g(t_{0})=0$ for some $t_{0}\geq 0$, then the existence of vanishing angles is guaranteed. For the research of vanishing calibration and the analysis for the solutions of $g(t)$, one can refer to \cite[Sects. 2 and 3]{La91}. The function $p(t)$ has a lower bound 
$$
F(\alpha, t, k-1):=\left(1-t \alpha \sqrt{\frac{k-2}{k-1}}\right)\left(1+t \alpha \sqrt{\frac{1}{(k-1)(k-2)}}\right)^{k-2}
$$
where $\alpha:={\rm sup}_{v}||A_{v}||$ and $F(\alpha, t, k-1)$ is a nonincreasing function of $k-1$, its limit as $k \rightarrow \infty$ is $(1-\alpha t)e^{\alpha t}$.

\medskip\noindent

Now, let $\theta_{1}(k,\alpha)$ and $\theta_{2}(k,\alpha)$ be the corresponding two estimated vanishing angles for the equation \eqref{va1} in equality form by replacing $p(t)$ with $F(\alpha, t, k-1)$ and $(1-\alpha t)e^{\alpha t}$ if they exist. Also from \eqref{new ODE}, one can easily see that the existence of $\theta_2$, $\theta_{1}$ implies the existence of $\theta_{0}$. Then
\begin{equation}\label{mid}
    \theta_0 \leq \theta_1(k, \alpha) \leq \theta_2(k, \alpha),
\end{equation}
and Lawlor made a "Table" (cf. \cite[Sect. 1.4]{La91}) of estimated vanishing angles by numerical analysis, including $\theta_1$ for $\operatorname{dim} C=\{3, \cdots, 11\}, \theta_2$ for $\operatorname{dim} C=12$. If ${\rm dim} \ C=k>12$, the following formula (\cite[Prop. 1.4.2]{La91}) for estimating vanishing angles was given by Lawlor
\begin{equation}\label{bigd}
  \tan \left(\theta_2(k, \alpha)\right)<\frac{12}{k} \tan \left(\theta_2\left(12, \frac{12}{k} \alpha\right)\right).
\end{equation}

\medskip\noindent

\textbf{Remark 3.6.} \textit{Based on Lawlor's method, many higher codimensional area-minimizing cones were found, as in \cite{P93}, \cite{Ke94}, \cite{Ka02}, \cite{XYZ18}, \cite{TZ20}, \cite{JC22}, \cite{JCX22}, \cite{CJL25}, \cite{Z25}, etc.}

\medskip\noindent
\medskip\noindent

\subsection{The second fundamental forms and the vanishing angles}

\medskip\noindent

Denote the second fundamental forms of  $\Sigma \hookrightarrow \mathbb{S}^{n} \times \mathbb{S}^{n}$ by $B$, denote the second fundamental forms of $\Sigma \hookrightarrow \mathbb{S}^{2n+1}(\sqrt{2})$ (and $\frac{1}{\sqrt{2}}\Sigma \hookrightarrow \mathbb{S}^{2n+1}(1)$) by $\widetilde{B}$. Denoted the two slices of $\Sigma$ by $\Sigma_{1}: \Sigma \cap\left(S^{n} \times\{y\}\right) \hookrightarrow S^{n} \times\{y\}$ and $\Sigma_{2}: \left(\{x\}\times S^{n}\right) \cap \Sigma \hookrightarrow\{x\} \times S^{n}$. We let $\bar\nabla^{1}$ and $\bar\nabla^{2}$ denote the Riemannian connections in $S^{n} \times\{y\}$ and $\{x\} \times S^{n}$. Recall the following canonical frame \eqref{F} along $\Sigma$:
\begin{equation}\label{cfe2}
    \begin{aligned}
      e_{0}=&\left(\frac{N_{1}}{|N_{1}|},0\right)=(\varepsilon_{0},0),e_{1}=(\varepsilon_{1},0),\ldots,e_{n-1}=( \varepsilon_{n-1},0),\\ e_{0}^{\prime}=&\left(0, \frac{N_{2}}{|N_{2}|}\right)=(0, \varepsilon_{0}^{\prime}),e_{1}^{\prime}=(0,\varepsilon_{1}^{\prime}), \ldots,e_{n-1}^{\prime}=(0,\varepsilon_{n-1}^{\prime}),\\ 
    \end{aligned}
\end{equation}
then $\Sigma_{1}$ (resp. $\Sigma_{2}$) has tangent frame $\{e_{i}\}_{i=1}^{n-1}$ (resp. $\{e^{\prime}_{i}\}_{i=1}^{n-1}$), $\Sigma$ has tangent frame  $\{e_{i}, e^{\prime}_{j}, T=\frac{e_{0}-e_{0}^{\prime}}{\sqrt{2}}\}$ at $(x,y) \in \Sigma$. We let $B^1$ (resp. $B^2$) denote the induced Riemannian connection and second fundamental form of $\Sigma_{1}$ with respect to $\varepsilon_{0}$ (resp. $\Sigma_{2}$ with respect to $\varepsilon_{0}^{\prime}$). 

\medskip\noindent

In this section, we compute the maximum value of the square norm for the shape operator of $\frac{1}{\sqrt{2}}\Sigma \subset  \mathbb{S}^{2n+1}(1)$ $$\alpha^2:=\sup _{\nu \in N \frac{\Sigma}{\sqrt{2}} \ \cap \  T \mathbb{S}^{2 n+1}(1)}\left\|\widetilde{B}^{\nu}\right\|^2,$$ then estimate the vanishing angles.

\medskip\noindent

Firstly, following Theorem 2.3, note that $m-1=(n-1)F$ on $\Sigma$,  one can get $||B||^2$. 

\medskip\noindent

\textbf{Proposition 3.7.} \textit{The length of the second fundamental forms of $\Sigma$ in $\mathbb{S}^n \times \mathbb{S}^n$ equals to $||B||^2=3(n-1)$.}

\medskip\noindent

Thus $\frac{1}{\sqrt{2}}\Sigma \hookrightarrow \mathbb{S}^{n} \left( \frac{1}{\sqrt{2}}\right) \times \mathbb{S}^{n} \left( \frac{1}{\sqrt{2}}\right)$ has $||B||^2=6(n-1)$. In the following, we show that $\alpha^2=6(n-1)$. 
Denote the two normal vector of $\Sigma$ in $\mathbb{S}^{2n+1}(\sqrt{2})$ as: $\nu_{1}=\frac{(x,-y)}{\sqrt{2}}$, $\nu_{2}=N=\frac{e_{0}+e_{0}^{\prime}}{\sqrt{2}}$.

\medskip\noindent

Firstly, the following Lemma is crucial in our analysis.

\medskip\noindent

\textbf{Lemma 3.8.} \textit{The slices $\Sigma_1$ and $\Sigma_2$ are also minimal isoparametric hypersurfaces in the spheres.}

\medskip\noindent

\textbf{Proof:} Isoparametric properties can be directly verified by the computations in the proof of Theorem 3.2; we now show that they are also minimal hypersurfaces. Let $c:=\frac{m-1}{n-1}$, that is equivalent to saying, the zero locus of the homogeneous polynomial 
$$
\widetilde{F}(x):=\langle x, y\rangle^2+\sum_{q=1}^{m-1}\langle A_q x, y\rangle^2-c|x|^2=0
$$ 
is a minimal cone for some fixed $y \in \mathbb{S}^n$. The minimality is equivalent to (cf. \cite{Hs67} or Lemma 2.2)
$$
\langle D|D\widetilde{F}|^2, D\widetilde{F}\rangle = 2|D\widetilde{F}|^2 \Delta^{E}\widetilde{F} \ \ ({\rm mod} \ \widetilde{F})
$$
where $D, \Delta^{E}$ are the Euclidean gradient and Laplacian operators. The computations are simple, and we omit the details here. It implies that ${\rm trace} B^1=0$. Similarly ${\rm trace} B^2=0$. $\Box$

\medskip\noindent

Now, we compute $\widetilde{B}^{\nu_2}$ and $\widetilde{B}^{\nu_1}$.

\medskip\noindent

By the decomposition of the Riemannian connections for product manifolds $\mathbb{S}^n \times \mathbb{S}^n$ (cf. \cite[p.139]{Car92}) or just consider choosing curves that move in the slices for every time of computations, one has (see Sect. 1.3 for the notations)
$$
\widetilde{B}^{\nu_2}\left(e_i, e_j\right)  = \langle \widetilde{\nabla}_{e_i} e_j, \nu_{2} \rangle =\langle \bar{\nabla}_{e_i} e_j, \nu_{2} \rangle= \left \langle \left(\bar\nabla^{1}_{\varepsilon_i} \varepsilon_j, 0\right),  \frac{1}{\sqrt{2}}e_{0}  \right \rangle =\frac{1}{\sqrt{2}} B^1\left(\varepsilon_i, \varepsilon_j\right)
$$
and  
$$
\widetilde{B}^{\nu_2}\left(e_i, T\right) =\langle \widetilde{\nabla}_{e_i} T, \nu_{2} \rangle =\langle \bar{\nabla}_{e_i} T, \nu_{2} \rangle =\left\langle\bar \nabla_{e_i}\left(\frac{e_0-e_0^{\prime}}{\sqrt{2}}\right), \frac{e_0+e_0^{\prime}}{\sqrt{2}}\right\rangle= \left\langle\bar \nabla_{e_i}e_{0},e_{0}\right \rangle =0
$$
similarly one gets $\widetilde{B}^{\nu_2}\left(e_i^{\prime}, e_j^{\prime}\right)=\frac{1}{\sqrt{2}} B^1\left(\varepsilon_i^{\prime}, \varepsilon_j^{\prime}\right)$, $\widetilde{B}^{\nu_2}(e_{i}^{\prime},T)=0$. And, those zero terms can also be obtained following $\bar\nabla C=-2A(T)$ \cite[Lemma 1(1)]{Ur19} since $C\equiv0$ for our examples here.

\medskip\noindent

By using the isoparametric properties (cf. \cite[Theorem 3.5]{CR15}) $\bar{\nabla}_{e_{0}}e_{0}=\bar{\nabla}_{e_{0}^\prime}e_{0}^\prime=0$, one gets 
$\widetilde{B}^{\nu_2}(T,T)=0$. The extra terms $\widetilde{B}^{\nu_2}\left(e_i, e_j^{\prime}\right)$ can just be denoted  by $A_{ij}$. 

\medskip\noindent

Thus
\begin{equation}\label{crt}
\widetilde{B}^{\nu_2}=\left(\begin{array}{ccc}
\frac{1}{\sqrt{2}} B^1 & A & 0 \\
A^{T} & \frac{1}{\sqrt{2}} B^2 & 0 \\
0 & 0 & 0
\end{array}\right). 
\end{equation}

\medskip\noindent

In the direction $\nu_{1}=\frac{(x,-y)}{\sqrt{2}}$, we calculate 
\begin{equation}\label{3ndb}
\begin{aligned}
\widetilde{B}^{\nu_1}\left(e_i, e_j\right) & =\left\langle \widetilde{\nabla}_{e_i} e_j, \nu_1 \right\rangle=  \left\langle D_{e_i} e_j, \nu_1 \right\rangle=-\left\langle D_{e_i} \nu_1, e_j \right\rangle=-\frac{1}{\sqrt{2}}\delta_{ij},   \\
\widetilde{B}^{\nu_1}\left(e_i^{\prime}, e_j^{\prime}\right) & =\frac{1}{\sqrt{2}}\delta_{ij},  \ \ \  \ \widetilde{B}^{\nu_1}\left(e_i, e_j^{\prime}\right)= 0, \\
\widetilde{B}^{\nu_1}\left(e_i, T\right) & = -\left\langle D_{e_i} \nu_1 , \frac{e_0-e_0^{\prime}}{\sqrt{2}} \right\rangle=0, \\
\widetilde{B}^{\nu_1}\left(e_i^{\prime},  T\right) & =0, \\
\widetilde{B}^{\nu_1}\left(T, T\right) & =\left \langle D_{\frac{e_0-e_0^{\prime}}{\sqrt{2}}}\left(\frac{e_0-e_0^{\prime}}{\sqrt{2}}\right), \nu_1 \right\rangle=0.
\end{aligned}
\end{equation}
Then,
\begin{equation}
\widetilde{B}^{\nu_1}=\left(\begin{array}{ccc}
-\frac{1}{\sqrt{2}} I & 0 & 0 \\
0 & \frac{1}{\sqrt{2}} I & 0 \\
0 & 0 & 0
\end{array}\right).
\end{equation}

\medskip\noindent

The minimalities of $\Sigma$ and its two slices implies that ${\rm trace} \widetilde{B}^{\nu_1} \cdot \widetilde{B}^{\nu_2}=0$. Then, for a given unit normal $\nu=\cos \beta \nu_1+ \sin \beta \nu_2$, we have $|\widetilde{B}^{\nu}|^2=\cos^2 \beta |\widetilde{B}^{\nu_1}|^2+\sin^2 \beta |\widetilde{B}^{\nu_2}|^2$. Now, $|\widetilde{B}^{\nu_1}|^2=n-1$ and $|\widetilde{B}^{\nu_2}|^2=3(n-1)$ from Proposition 3.7[]
 thus $\sup _\nu\left|\widetilde{B}^{\nu}\right|^2=3(n-1)$ and finally
$$
\alpha^2=\sup _{\nu \in N \frac{M}{\sqrt{2}} \ \cap \  T \mathbb{S}^{2 n+1}(1)}\left\|\widetilde{B}^{\nu}\right\|^2=6(n-1).
$$

\medskip\noindent

\textbf{Proposition 3.9.} \textit{$\alpha^2=6(n-1)$.}

\medskip\noindent

Now, $({\rm dim}\ C, \alpha^2)=(2n,6n-6)$, by checking Lawlor's table (\cite[p.20-21]{La91}), one finds that the vanishing angles don't exist for the case $3\leq 2n \leq 12$. 

\medskip\noindent

For the cases $2n > 12$, we apply Lawlor's formula \eqref{bigd} (\cite[Prop. 1.4.2]{La91}). The function 
$$
f(n):=\frac{12}{k}\alpha=\frac{6\sqrt{6}\sqrt{n-1}}{n}
$$
is decreasing in $n$, and $f(10)^2=19.44$, $f(11)^2\approx17.85$. By checking Lawlor's table for  ${\rm dim}\ C=12$, one can find that if $n\geq11$, then the vanishing angles exist and 
\begin{equation}\label{bigd2}
  \tan \left(\theta_2(2n, \sqrt{6}\sqrt{n-1})\right)<\frac{12}{2n} \tan 11.23^{\circ}<\frac{1.2}{n}.
\end{equation}

\medskip\noindent
\medskip\noindent

\subsection{Normal radius: searching for the shortest normal geodesics}

\medskip\noindent

Now, we compute the normal radius of the cone $\mathcal{C}$ over $\Sigma$. Firstly, 
$$\mathcal{C}=\{(x,y) \in  \mathbb{R}^{n+1} \times \mathbb{R}^{n+1}: |x|^2-|y|^2=0 \ \ {\rm and} \ \ \langle x, y\rangle^2+\sum_p B_p^2=c|x|^2 |y|^2 \},$$
where $c=\frac{m-1}{n-1}$, $B_p=\langle A_p x, y \rangle=\left\langle A_p^{\top} y, x\right\rangle$.

\medskip\noindent

Along $\Sigma$, denote $\overrightarrow{n_2}=(x,-y)$, it is of length $\sqrt{2}$. Denote 
$$
Y=\left( \langle x, y\rangle y+\sum_p B_p A_p^{\top} y-c x, \langle x, y\rangle x+\sum_p B_p A_p x-c y\right),
$$
then along $\Sigma$, $|Y|^2=2c(1-c)$, let $c_1:=\frac{1}{\sqrt{c(1-c)}}$, we set $\vec{n}_1=c_1 Y$. Then a normal geodesic of  $\Sigma$ located in $\mathbb{S}^{2n+1}(\sqrt{2})$ which starts from $(x,y)$ in the direction $v=\cos \alpha \vec{n}_1+\sin \alpha \vec{n}_2$  can be given as
$$
\gamma_v(\theta)=\cos \theta(x, y)+\sin \theta \cos \alpha(x,-y)+c_1\sin \theta \sin \alpha \ Y .
$$

\medskip\noindent

Denote
\begin{equation}
\left\{\begin{array}{l}
t=c_1 \sin \theta \sin \alpha \\
r=\cos \theta+\sin \theta \cos \alpha-c t \\
s=\cos \theta-\sin \theta \cos \alpha-c t
\end{array}\right.
\end{equation}
then 
$$
 \gamma_v(\theta)=(\tilde{x}, \tilde{y}) =
 \left(r x+t\langle x, y\rangle y+t \sum_p B_p A_p^{\top} y, 
\ \ s y+t\langle x, y\rangle x+t \sum_p B_p A_p x \right).
$$

We assume $\gamma_v(\theta)$ intersects $M$ with another point, then 
\begin{equation}\label{fie}
|\tilde{x}|^2=|\tilde{y}|^2=1
\end{equation}
and
\begin{equation}\label{see}
\langle\tilde{x}, \tilde{y} \rangle^2+\sum_q\left\langle A_q \tilde{x}, \tilde{y}\right\rangle^2=c
\end{equation}

\medskip\noindent

By using the properties of the skew-symmetric representations of the Clifford algebra in \eqref{skc}, after a series of complex computations, \eqref{fie} is equivalent to \begin{equation}
\begin{aligned}
\left\{\begin{array}{l}
r^2+t^2 c+2 r t c=1 \\
s^2+t^2 c+2 s t c=1
\end{array}\right.
\end{aligned}
\end{equation}
which implies that $2(r-s) \cos \theta=0$. Since $\theta \neq 0$. Then, it happens that, case (i): $r \neq s, \Rightarrow \cos \theta=0, \theta=\frac{\pi}{2}$; case (ii): $\quad r=s$. i.e. $\sin \theta \cos \alpha=0$, direct computations show that (1) is naturally satisfied.
Then, we now consider the case $\cos \alpha=0$, i.e. $v= \pm \overrightarrow{n_2}$. In this case, by repeatedly using the properties of the skew-symmetric representations of the Clifford algebra in \eqref{skc}, the equation \eqref{see} is reduced to
$$\cos ^2 \theta+2 \tan \theta(1-c)-\left(c-c^2\right) t^2=\cos 2 \theta \pm \sqrt{\frac{1-c}{c}} \sin 2 \theta=\pm 1.$$

\medskip\noindent

By setting $\varphi \in (0, \frac{\pi}{2})$, $\tan \varphi=\sqrt{\frac{1-c}{c}}$, we find that $\cos(2\theta \pm \varphi)=\cos(k\pi \pm \varphi)$ where two $\pm$ are independent. It leads to $\theta=\frac{k\pi}{2}-\varphi$ or $\theta=\varphi+\frac{k\pi}{2}$ or $\theta=\frac{k\pi}{2}$ for $k \in \mathbb{Z}$. By searching for the minimum value of $\theta \in (0,\pi)$, we find that the normal radius is 
$$
N= \min \left\{\operatorname{arctan} \sqrt{\frac{c}{1-c}}, \operatorname{arctan} \sqrt{\frac{1-c}{c}}\right\}
$$
where $c=\frac{m-1}{n-1}$.

\medskip\noindent

\textbf{Remark 3.10.} \textit{Generally, we only need to check a weaker condition (cf. \cite[Theorem 1.2.1]{La91}) that the normal wedges of angular radii equal to the vanishing angles are disjoint. For the isoparametric cone $\mathcal{C}$ over an isoparametric hypersurface $\Sigma$ in the sphere, two normal wedges whose radii are smaller than the focal radius will never intersect, since the parallel hypersurfaces of $\Sigma$ are also embedded isoparametric hypersurfaces if the distance is less than the focal distance. Also, the positivity of $p(t)=\mathrm{inf}_{v}{\rm det}(I-t A_v)$ implies that the resulting vanishing angles are always less than the focal radius. Thus, we can only consider the existence of vanishing angles for isoparametric cones. However, what we consider here are isoparametric hypersurfaces in $\mathbb{S}^n \times \mathbb{S}^n$, which are then embedded into $\mathbb{S}^{2n+1}(\sqrt{2})$ as codimension-two submanifolds. For this reason, the computations in this section are necessary.} 

\medskip\noindent
\medskip\noindent

\subsection{The proof of Theorem B}

\medskip\noindent

\textbf{Proof of Theorem B}: we are concerned with the cases $m\geq 2$, for $m=1$, these cones are complex varieties, naturally area-minimizing. As in Sect. 2.2, the restrictions for the multiplicity of FKM isoparametric hypersurfaces in the sphere $\mathbb{S}^{2n+1}$ satisfy: $n\geq1, n-m \geq 1$, one can find that the first family of FKM isoparametric hypersurfaces in $\mathbb{S}^{n} \times \mathbb{S}^{n}$ appear in $\mathbb{S}^{3} \times \mathbb{S}^{3}$. Our requirement for the existence of vanishing angles is that $n= k\delta(m)-1 \geq 11$, and we also have $k\delta(m)-1-m\geq 1$, i.e. $k\delta(m) \geq m+2$. The normal radius $N=\arctan\sqrt{\frac{c}{1-c}}$ if $c=\frac{m-1}{n-1}\leq \frac{1}{2}$ which is equivalent to $k\delta(m) \geq 2m$, oterwise $N=\arctan \sqrt{\frac{1-c}{c}}$. 

\medskip\noindent

\textbf{Case (i): $k=1$.}  We need $m\geq 2, \delta(m)\geq m+2, \delta(m)\geq 12$, and decide whether $\delta(m)\geq 2m$. By checking table \eqref{table}, we need $m\geq 9$. 

\medskip\noindent

If $m=9$, then $n=15$. The normal radius now is $\arctan\sqrt{\frac{1-c}{c}}=\arctan{\frac{\sqrt3}{2}}$, which is bigger than our upper bound for the two times of estimated vanishing angles: $2\arctan \frac{1.2}{n}= 2\arctan \frac{1.2}{15}$.

\medskip\noindent

If $m\geq 10$, then $n=\delta(m)-1\geq 31$.
the normal radius is 
$$
\arctan\sqrt{\frac{c}{1-c}}=\arctan\sqrt{\frac{m-1}{n-m}}\geq \arctan\sqrt{\frac{9}{n-10}},
$$
then if 
$$
\arctan\sqrt{\frac{9}{n-10}}>2\arctan \frac{1.2}{n}
$$
which is equivalent to 
$$
n^2-1.44-2.4n\sqrt{\frac{n-10}{9}}>0
$$
two times of vanishing angles are all less than the associated normal radius, and it is naturally satisfied for all $n\geq 11$. Together with the case $m=9$, we lead to $n\geq 15$.

\medskip\noindent

\textbf{Case (ii): $k\geq 2$.} Then naturally it satisfy $k\delta(m) \geq 2m$, the normal radius is 
$$
N=\arctan\sqrt{\frac{c}{1-c}}= \arctan\sqrt{\frac{m-1}{n-m}}
$$
where $n=k\delta(m)-1$. The conditions $n\geq 11, n\geq m+1$ can  satisfied for all $m\geq 2$ with suitable given $k$. Thus if 
\begin{equation}\label{k2}
    2 \arctan \frac{11}{n}<\arctan\sqrt{\frac{1}{n-2}}
\end{equation}
two times of vanishing angles are all less than the associated normal radius since $\frac{m-1}{n-m}$ is increasing in $m$. \eqref{k2} is equivalent to 
$$
n^2-1.44-2.4n\sqrt{n-2}>0,
$$
and it is also satisfied for all $n\geq 11$. And $n=11$ can be attained for $m=3,k=3$ or $m=4,k=3$.
$\Box$

\medskip\noindent
\medskip\noindent

\section{Area-minimizing cones of minimal product types}

\medskip\noindent

We let $f_i: \Sigma_i\hookrightarrow S^{a_i}(1)$ be a family of minimal immersions with ${\rm dim} \Sigma_{i}=k_{i}$ for $1 \leq i \leq m$ and $k_1\leq k_{2}\leq \cdots \leq k_{m}$, ${\rm dim} \Sigma =k-1=\sum_{i=1}^{m}k_{i}$. The product immersion $f$ is
$$
\begin{aligned}
f: \Sigma=\Sigma_1 \times \cdots \times \Sigma_m & \rightarrow S^{a_1+\cdots+a_m+m-1}(1) \\
\left(x_1, \ldots, x_m\right) & \mapsto\left(\lambda_1 f_1\left(x_1\right), \ldots, \lambda_m f_m\left(x_m\right)\right),
\end{aligned}
$$
where $\sum_{i=1}^m \lambda_i^2=1$. Then, $f$ is minimal if and only if (cf. \cite{CH18}, \cite{TZ20} and \cite{JCX22}):

$$
\lambda_i=\sqrt{\frac{k_i}{\sum_{i=1}^m k_i}}.
$$

\medskip\noindent

We call the cones $\mathcal{C}(\Sigma)$ minimal product cones. And area-minimizing minimal product cones provide new topological types for area-minimizing cones that can't be split into products of area-minimizing cones. To study the area-minimization of minimal product cones via Lawlor's curvature criterion, the following formulas for second fundamental forms and normal radius are needed. We remark here, using these formulas, Zhang recently (\cite[Theorem 1.1]{Z25}) observed an important property for minimal product cones, that the minimal product cone over the $n$-copies of an arbitrary embedded minimal submanifold in the unit sphere can be area-minimizing if $n$ is sufficiently large, this is an amazing fact. The main reason is that, by checking Prop 4.1 and Prop 4.2, one can find that the convergence rate of the normal radius is of order $O\left(\frac{1}{\sqrt{n}}\right)$, while that of the vanishing angle is $O\left(\frac{1}{n}\right)$. And, one can see \cite{Z25} for more results on area-minimizing cones of minimal product type.

\medskip\noindent

\textbf{Proposition 4.1 (\cite[Equ. (32)]{TZ20} and \cite[Theorem 2.2]{JCX22}).} \textit{The upper bound for the norms of shape operator of $\Sigma$ is given by $\alpha^2=\operatorname{dim} \Sigma \cdot \max \left\{1, \frac{\alpha_1^2}{k_1}, \ldots, \frac{\alpha_m^2}{k_m}\right\}$, where $\alpha_i^2$ are the upper bound for the norm of shape operator of $\Sigma_{i}$.}

\medskip\noindent

We now let all $\Sigma_{i}$ be the minimal FKM isoparametric hypersurfaces in $\mathbb{S}^{n_{i}} \left( \frac{\sqrt{2}}{2} \right) \times \mathbb{S}^{n_{i}} \left( \frac{\sqrt{2}}{2} \right)\subset \mathbb{S}^{2n_{i}+1}(1) $ with $n_{i}\geq 3$. Then as codimensional 2 submanifolds, ${\rm dim} \Sigma_{i}=k_{i}=2n_{i}-1$ and 
 $\alpha_i^2=6n_{i}-6=3k_i-3$ as in Sect. 5. Then we find $\alpha^2= \operatorname{dim} \Sigma \left(3-\frac{3}{k_{m}}\right)=3(k-1)(1-\frac{1}{k_{m}})$ for $\mathcal{C}(\Sigma)$.  

 \medskip\noindent
 
Then by using Lawlor's formula \eqref{bigd} and checking the last column in Lawlor's table (\cite[p. 21]{La91}), we find that if 
$$
\left(\frac{12\sqrt{3(k-1)(1-\frac{1}{k_{m}})}}{k}\right)^2<19
$$
then the estimated vanishing angles exist, by letting $k_{m}$ tend to $k-1$, we find that if $k\geq 21$, then the estimated vanishing angles $\theta_2$ exist and satisfy   
 $$
 \theta_{2} \leq {\rm arctan} \frac{5}{2k}.
 $$

\medskip\noindent

The following formula for normal radius was obtained independently in \cite[Theorem 3.3]{JCX22} and \cite[Corollary 3.5]{Z25}.

\medskip\noindent

\textbf{Proposition 4.2 (\cite[Theorem 3.3]{JCX22} and \cite[Corollary 3.5]{Z25}).} \textit{Given a minimal product $\Sigma=\Sigma_1 \times \cdots \times \Sigma_m(m \geq 2)$ over submanifolds which are not totally geodesic spheres. Then the normal radius $N$ of $\Sigma$ (or its cone) at $x=\left(\lambda_1 x_1, \cdots, \lambda_m x_m\right)$ satisfy
$$
\operatorname{dim} \Sigma (1- \cos N)=\min _{1 \leq i \leq m}
\operatorname{dim} \Sigma_i\left(1-\cos N_i\right)
$$
where $\Sigma_i$ are the normal radius of $\Sigma_i$ at $x$}.

\medskip\noindent

Since $$
{\rm cos} \ N_{i}= \max \left\{\sqrt{1-c_{i}}, \sqrt{c_{i}}\right\},
$$
where $c_{i}=\frac{m_{i}-1}{n_{i}-1}\geq \frac{1}{n_{i}-1}$, then ${\rm cos} \ N_{i} \leq \sqrt{\frac{n_{i}-2}{n_{i}-1}}$. And finally, we estimate that
$$
\operatorname{dim} \Sigma (1- \cos N) \geq \min _{1 \leq i \leq m}  (2n_{i}-1) \left(1- \sqrt{\frac{n_{i}-2}{n_{i}-1}}\right) \geq 1
$$
since the right hand is decreasing in $n_{i}\geq 3$ and tends to $1$ as $n_{i}$ tends to infinity, then the normal radius has a lower bound 
$$
N \geq {\rm arccos}\left( 1-\frac{1}{k-1}\right).
$$

It is easy to check that when $k \geq 21$, 
$$
{\rm arccos}\left( 1-\frac{1}{k-1}\right) \geq 2{\rm arctan} \frac{5}{2k}
$$
is always satisfied by considering the following inequalities: $2{\rm arctan} \frac{5}{2k}<\frac{5}{k}, 1-\cos \frac{5}{k}<\frac{25}{2k^2}$. Then, two times of vanishing angles are less than normal radius. Thus, we obtain that 

\medskip\noindent

\textbf{Theorem 4.3.} \textit{The minimal product cone $\mathcal{C}$ over FKM isoparametric hypersurfaces is area-minimizing if ${\dim C} \geq 21$.}

\medskip\noindent

\textbf{Remark 4.4.} \textit{The numbers 16, 12 and 21 are uniform lower bounds. However, they could be further improved if one performs more numerical computations than Lawlor did, as his table (\cite[p.20-21]{La91}) only includes values for $3\leq{\dim \mathcal{C}} \leq 12$.}

\medskip\noindent
\medskip\noindent

\section{Comments and questions}

\medskip\noindent

The discovery of the new isoparametric foliation structure is precious and may have many potential applications, see \cite{GQTY25} for an extensive review on their applications. We here discuss a possible application in singular Riemannian foliations of spheres. In \cite{R14}, M. Radeschi cosntruct new indecomposable singular Riemannian foliations on round spheres, by the map 
$$
\begin{aligned}
\pi_C: \mathbb{S}^{2 n+1} & \longrightarrow \mathbb{R}^{m+1} \\
x & \longmapsto\left(\left\langle P_0 x, x\right\rangle, \ldots\left\langle P_m x, x\right\rangle\right) .
\end{aligned} 
$$
where $\{P_{i}\}$ are a family of Clifford system. It can be checked that the following map is an isoparametric map in the spirit of Wang (\cite{W86}, \cite{W94})
$$
\begin{aligned}
\pi: \mathbb{S}^{2 n+1} & \longrightarrow \mathbb{R}^{2} \\
(x,y) & \longmapsto\left( |x|^2-|y|^2, F(x,y) \right)
\end{aligned}
$$
where $F(x,y)$ is the splitting form of the FKM isoparametric polynomial, it is the square sum of the components in Radeschi's definition up to a constant. The regular fibers are isoparametric hypersurfaces in  $\mathbb{S}^{n}(a) \times \mathbb{S}^{n}(b) \subset \mathbb{S}^{2 n+1}(1)$ which has codimensional $2$. However, it can be directly checked that, the Lie bracket of spherical gradients is not zero, thus $\pi$ is not an isoparametric map of Terng \cite{T85}, as pointed in \cite{T10}, the preimages of Terng's isoparametric map gives singular Riemannian foliations, thus we are interesting in that whether Wang's isoparametric map $\pi$ gives singular Riemannian foliations in $\mathbb{S}^{2 n+1}$. 

\medskip\noindent

For other potential applications, a study on solutions for Yamabe-type problems in $\mathbb{S}^{n} \times \mathbb{S}^{n}$ realted to Qian-Tang and Urbano's isoparametric hypersurfaces is given in \cite{HJ25}. Searching for closed embedded self-shrinkers of mean curvature flow by reduction method of isoparametric maps, as in \cite{R23}, \cite{MM24} etc.

\medskip\noindent

The inhomogeneous properties of FKM isoparametric hypersurfaces in $\mathbb{S}^{2n+1}$ was studied very well as in \cite{FKM81} (also see \cite{CR15}), the isometric group of $\mathbb{S}^n \times \mathbb{S}^n$ is the subgroup of the orthogonal group $O(2n+2)$ given by
$$
\left\{\left(\begin{array}{cc}
A & 0 \\
0 & B
\end{array}\right),\left(\begin{array}{cc}
0 & A \\
B & 0
\end{array}\right) / A, B \in \mathrm{O}(n+1)\right\}.
$$
we thus pose the following question

\medskip\noindent

\textbf{Question 1.} \textit{Are these non-homogeneous FKM isoparametric hypersurfaces in $\mathbb{S}^{2n+1}$ necessarily non-homogeneous when restricting on $\mathbb{S}^{n} \times \mathbb{S}^{n}$?}

\medskip\noindent

So far, our discussions in this paper concern isoparametric hypersurfaces. The geometries of focal varieties are clearly important in isoparametric theory. By the well-known fact (cf. \cite[Theorem 1.1]{GT14}), both focal varieties are minimal submanifolds. $F^{-1}(0)$ is the same as FKM examples which was well-studied before, however, $F^{-1}(1)$ is the intersection of focal varieties of FKM isoparametric hypersurfaces in spheres with $\mathbb{S}^{n}(1) \times \mathbb{S}^{n}(1)$. The cones over $F^{-1}(0)$ was well-studied in \cite{TZ20}. Thus, we prefer to ask the following question:

\medskip\noindent

\textbf{Question 2.} \textit{Are the focal submanifolds $F^{-1}(1)$ minimal in $\mathbb{S}^{2n+1}(\sqrt{2})$? If so, what can be said about the area-minimizations of their cones? Can we define similar product angle function $C$ for high codimensional submanifolds of $\mathbb{S}^{n}(1) \times \mathbb{S}^{n}(1)$?}

\medskip\noindent

\textbf{Wang's isoparametric maps and reduction method.} There are some low-dimensional families (cf. \eqref{table}) in Theorem B which remain unsolved, due to the limitations of Lawlor's method in low-dimensional cone cases. In fact, it cannot deal with the cones over Qian-Tang and Urbano's examples in low-dimensional cases.
In contrast to Wang's work in \cite{W94}, the classification of area-minimizing isoparametric cones depends on the reduction method for rank $2$ isoparametric map (also see \cite[Chapter 4]{ER93}), however, what Wang considered are hypersurfaces in the sphere, whereas ours are codimension-two submanifolds of the sphere, thus, we need to find some rank 3 isoparametric map in the spheres, such that the preimage of geodesics in the $3$-dimensional reduction spaces are just our examples. This is not easy. Thus, we pose the question

\medskip\noindent

\textbf{Question 3.} \textit{How about their area-minimizations of the cones over left minimal isoparametric hypersurfaces of $\mathbb{S}^{n}(1) \times \mathbb{S}^{n}(1) \subset \mathbb{S}^{2n+1}(\sqrt{2})$ in Theorem B?}

\medskip\noindent

Our area-minimizing cones are constructed by the intersection of FKM isoparametric cones and Simons cones, as in Theorem B. We are also interested in finding new intersection varieties that are area-minimizing. For given intersection varieties, the minimality is not hard to verify, as in \cite{Hs67}. However, the classification problem is hard, and is still open for minimal cubic cones (see \cite[Remark 1.2]{CJX24} for an introduction).

\medskip\noindent

\textbf{Question 4.} \textit{Are there other area-minimizing intersection varieties given by two other types of isoparametric polynomials?}

\medskip\noindent
\medskip\noindent

\section{Appendix: regular area-minimizing cones of low codimension $c\leq 2$}

\medskip\noindent

Firstly, by Simons' famous result \cite{S68}, area-minimizing cones of codimension $c=1$ only exists in $\mathbb{R}^{l+1}$ for $l+1\geq 8$. We now review some famous examples except for Simons cones $\mathcal{C}_{n},n\geq3$. By a classical reduction technique \cite{HsL71} for homogeneous minimal hypersurfaces, together with a calibration method in the $2$-dimensional orbit spaces, Lawson \cite{La72} proved many homogeneous minimal hypersurfaces in the classification (cf. \cite{HsL71}) gave rise to area-minimizing cones, which including $\mathcal{C}_{n}(n\geq 3)$ and $C(\mathbb{S}^2 \times \mathbb{S}^4)$, etc. Subsequently, P. Simões \cite{Si74} proved that  $C(\mathbb{S}^1 \times \mathbb{S}^5)$ is stable but not area-minimizing. A complete classification for homogeneous area-minimizing hypercones was given in \cite[Sect. 5.1]{La91}. Y. Zhang \cite{Z16} improved
Lawson’s original calibration method and provided a new proof for this complete classification.

\medskip\noindent

In fact, up to now, all known regular area-minimizing hypercones are isoparametric cones (cf. \cite[Sect. 6]{S89}): cones over isoparametric hypersurfaces in the sphere $\mathbb{S}^l$ that have constant principal curvatures. According to the elegant classification results in \cite{CCJ07} and \cite{Chi20}: all non-homogeneous examples arise from the FKM isoparametric hypersurfaces (cf. \cite{FKM81} and see Sect. 2.2) in $\mathbb{S}^{2n+1}$. However, area-minimizing isoparametric cones were classified earlier: based on the characteristic foliation constructed by Ferus and Karcher \cite{FK85}, an explicit list is given in \cite [Sect. 9]{W94}\footnote{Wang's proof include an key reduction method that generalized Hsiang-Lawson's reduction theorem for orbit cases. For orbits, the reduction method means that homogeneous minimality can be restricted to critical points of equivariant variations \cite[Theorem 1]{HsL71}. For isoparametric map $f$, Wang \cite{W86} proved that the minimality of $f$-invariant hypersurfaces is also equivalent to being a critical point of $f$-invariant variations, see also \cite[Chapter 4]{ER93}.}, see \cite{W94} and \cite{TZ20} for detailed discussions of area-minimizing cones associated to isoparametric foliations in spheres. As a noteworthy application, L. Simon \cite{S89} constructed many nonlinear entire solutions to the minimal surface equation on $\mathbb{R}^n \ (n\geq 8)$ via area-minimizing isoparametric cones.

\medskip\noindent

Searching for other types of regular area-minimizing hypercones remains a challenging task. There is an open problem asked by S. T. Yau \cite[Problem 104]{Yau82}: classify the topological type of the $7$-dimensional (necessarily regular) area-minimizing cones in $\mathbb{R}^8$. In \cite{IW15}, Ilmanen and White asked whether cones over minimal (non-totally geodesic) embedded $\mathbb{S}^{l-1} \subset \mathbb{S}^{l}$ constructed in \cite{Hs83a}, \cite{Hs83b}, and \cite{HS86} are area-minimizing (see also Tomter's examples in \cite{T87}). 

\medskip\noindent

A few classes of regular area-minimizing cones of codimension $c=2$, which are non-holomorphic, were found before, and both are homogeneous examples. In contrast to real varieties, holomorphic cones with specific degree may not be strictly stable, even they are area-minimizing (cf. \cite{DL25}). In 1988, via a delicate generalization of Lawson's calibration method in codimension-two cases, Cheng \cite{Cheng88} found one family (and two isolated examples) of homogeneous codimension two area-minimizing cones. Shortly after Cheng, Lawlor \cite{La91} provided a general criterion for regular area-minimizing cones, known as the curvature criterion (see Sect. 3.1). Among his many results, the following homogeneous codimension-two area-minimizing cones (cf. \cite[Theorem 5.1.1]{La91}) were found (see his introduction in \cite{La91} for more historical notes):
$$
\mathcal{C}(\mathbb{S}^2\times \mathbb{S}^2\times \mathbb{S}^2) \ \ {\rm and} \ \ \mathcal{C}(\mathbb{S}^p\times \mathbb{S}^q\times \mathbb{S}^r), p+q+r\geq 7.
$$

\medskip\noindent

In the lowest-dimensional case, there is an open question posed by Lawlor (cf. \cite{La91}): whether the cone over Veronese embedded $\mathbb{R}P^2 \subset \mathbb{S}^4$ is area-minimizing.  Meanwhile, the codimensional $3$ area-minimizing cones over Veronese embedded $\mathbb{C}P^2 \subset \mathbb{S}^7$ play key roles in Liu's generic singularity study (cf. \cite[Sect. 8]{L25b}). 

        \medskip\noindent

 \textbf{Acknowledgements} The author wishes to express his sincere gratitude to Teng Wang and Prof. Xiaowei Xu for bringing the examples in Remark 2.10 to his attention at a very early stage. Building upon this initial information, the author subsequently established these more general FKM isoparametric examples, motivated by the key observation that their examples admit a formulation in Clifford form. Particularly, the author thanks Prof. Xiaowei Xu for drawing his attention to Reference \cite{Ur19} and many inspiring discussions. The author also wishes to thank Huixin Tan, Prof. Yuquan Xie, and Prof. Wenjiao Yan for drawing his attention to Reference \cite{QT16}, and thank Dr. Zeke Yao for his discussions on homogeneous hypersurfaces. This work is supported by the NSFC (No.12301068, No. 12571057), the project of Stable Support for Youth Team in Basic Research Field, CAS (YSBR-001), and Xiaomi Young Scholar Fund.

\medskip\noindent




	\medskip\noindent
		\medskip\noindent

\begin{flushleft}
			\medskip\noindent
			\begin{tabbing}
				XXXXXXXXXXXXXXXXXXXXXXXXXX*\=\kill
				Hongbin Cui\\
				School of Mathematical Sciences, University of Science and Technology of China\\
				Wu Wen-Tsun Key Laboratory of Mathematics, USTC, Chinese Academy of Sciences\\
				96 Jinzhai Road, Hefei, 230026, Anhui Province, China\\
				
				E-mail: cuihongbin@ustc.edu.cn
				
			\end{tabbing}
		\end{flushleft}

  \end{document}